\documentclass{article}
\usepackage{amsmath}
\def \E{\hbox{ I\hskip -2pt E}}

\def \P{\hbox{I\hskip -2pt P}}
\def \N{\hbox{\it I\hskip -2pt N}}
\def \C{\hbox{\it I\hskip -5.5pt C}}
\def\Pa{{{\cal P}}}
\def\Ca{{{\cal C}}}
\def \R {\hbox{\it I\hskip -2pt R}}
\def\Z{\hbox{\it Z\hskip -3pt Z}}
\def\squarebox#1{\hbox to #1{\hfill\vbox to #1{\vfill}}}

\title{First order asymptotics of matrix integrals ; a rigorous 
approach towards the understanding of matrix models  }
\author{
ALICE GUIONNET\thanks{Ecole Normale Superieure de Lyon,
Unite de Mathematiques pures et appliquees,
UMR 5669,
46 Allee d'Italie, 
69364 Lyon Cedex 07, France},\ 
\,\,
}

\date{30/10/2002}

\newtheorem{prop}{Proposition}[section]
\newtheorem{lem}[prop]{Lemma}

\newtheorem{pr}[prop]{Property}
\newtheorem{cor}[prop]{Corollary}
\newtheorem{theo}[prop]{Theorem}

\newtheorem{disc}[prop]{Discussion}

\newenvironment{rem}{\refstepcounter{prop}
{\bf{Remark }}\theprop :}{}

\newcommand{\eqnsection}{
   \renewcommand{\theequation}{\thesection.\arabic{equation}}
   \makeatletter
   \csname @addtoreset\endcsname{equation}{section}
   \makeatother}

\def\nn{\noindent}

\def\xx{\vrule height 0.7em depth 0.2em width 0.5 em}

\def\a{\alpha}
\def\b{\beta}
\def\d{\delta}
\def\e{\epsilon}

\def\l{\lambda}

\def\s{\sigma}

\def\D{\Delta}

\def\L{\Lambda}

\def\nn{\noindent}
\def\part{\partial}

\def\ts{\times}

\def\ra{\rightarrow}

\def\lbc{\lbrace}
\def\rbc{\rbrace}

\def\tilde{\widetilde}

\def\Aa{{\cal A}}

\def\CC{{\cal C}}

\def\Fa{{\cal F}}

\def\Ha{{\cal H}}

\def\La{{\cal L}}
\def\Ma{{\cal M}}

\def\Oa{{\cal O}}
\def\PP{{M_1}}

\def\Sa{{\cal S}}

\def\CC{{\cal C}}

\def\m1{{{\cal M}_1^+(\R)}}

\def\mun{{\hat \mu ^N}}

\def\prf{{\bf Proof.}}

\newcommand{\beaa}{\begin{eqnarray*}}
\newcommand{\eeaa}{\end{eqnarray*}}
\newcommand{\bea}{\begin{eqnarray}}
\newcommand{\eea}{\end{eqnarray}}
\newcommand{\be}{\begin{equation}}
\newcommand{\ee}{\end{equation}}

\def\11{{\hbox{1\kern-.2em\hbox{I}}}}
\def\tr{{\hbox{\rm tr}}}
\def\trn{{\hbox{\rm tr}_N}}

\def\on2{{\frac{1}{2}}}
\def\boxplus{{\fbox{+}}}

\begin{document}
\eqnsection
\maketitle

\nn
{\bf Abstract :} We investigate the limit behaviour 
of the spectral measures of matrices following the Gibbs measure
for the Ising model 
on random graphs, Potts model on random graphs, matrices
 coupled in a chain
model or induced QCD model. For most of these models, we prove 
that the  spectral measures  converge almost surely and describe 
their limit via solutions to an Euler equation for
isentropic flow with negative pressure $p(\rho)=-3^{-1}\pi^2 \rho^3$.

\bigskip

\medskip

\nn

\medskip

\nn
{\it Mathematics Subject of Classification :}  60F10, 15A52, 35F20.

\section{Introduction}
It appears since the  work of 't Hooft
that  matrix integrals can be seen, via Feynman diagrams expansion, 
as generating functions 
for enumerating maps (or triangulated surfaces).
We refer here to the very nice survey of 
A. Zvonkin's  \cite{ZVON}. One matrix integrals 
are used to enumerate maps with a given genus 
and given  vertices degrees distribution
whereas several matrices integrals can be used 
to consider the case where the vertices can additionally be
coloured (i.e can take different states).

Matrix integrals are usually of the following form

$$Z_N(P)=\int e^{-N\tr(P(A_1^N,\cdots,A_d^N))}
d
A_1^N\cdots dA_d^N$$
with some polynomial function $P$ of
$d$-non-commutative variables and 
the Lebesgue measure $dA$ on some well chosen  ensemble
of $N\ts N$ matrices such as the set $\Ha_N$ (resp.
$\Sa_N$, resp. $\Sa ymp_N$)
of $N\ts N$ Hermitian  (resp. symmetric, resp. symplectic) 
 matrices.  One would like to
understand the full expansion of
$Z_N(P)$ in powers of $N$. 
For instance, in the case where the
matrices live on $\Ha_N$, the formal expansion 
linked with Feynamn diagrams is of the type

$${1\over N^2} \log Z_N(P)=\sum_{g\ge 0} {1\over N^{2g}} 
C_P(g)$$
where $C_P(g)$ enumerates some maps with genus $g$.
Such an expansion was proved to hold rigorously 
in the one matrix case 
 by K. McLaughlin and N. Ercolani in 2002.

A related issue is to understand the 
asymptotic behaviour of the corresponding Gibbs
measure

$$\mu^N_P(d
A_1^N\cdots dA_d^N)={1\over Z_N(P) }
 e^{-N\tr(P(A_1^N,\cdots,A_d^N))}
d
A_1^N\cdots dA_d^N.$$
More precisely, if for a $N\ts N$ matrix $A$, 
$(\l_1(A),\cdots,\l_N(A))$ denotes its eigenvalues 
and $\mun_A:=N^{-1}\sum_{i=1}^N\d_{\l_i(A)}$
its spectral measure, one would like to
understand the asymptotic 
behaviour of $(\mun_{A_1^N},\cdots,\mun_{A_d^N})$
under the Gibbs measure $\mu^N_P$ when
$N$ goes to infinity. Of course, this understanding is
intimately related with the
first order asymptotic
of the free energy $F_N(P)=N^{-2}\log Z_N(P)$.
In fact, the rigorous approach 
of the full expansion of matrix integrals 
when $d=1$ given by K. McLaughlin and N. Ercolani
is based on Riemann Hilbert problems techniques which
themselves require a precise understanding 
of such asymptotics of the spectral measures.

However, only very few 
 matrix integrals could be evaluated 
in the physics litterature, even on a non rigorous ground.
These cases corresponds in general to the
case where integration holds over Hermitian matrices.
Using orthogonal polynomial methods, Mehta \cite{mehta}
obtained the limiting free energy for the Ising model
on random graphs, corresponding to $d=1$ and $P(A,B)= P(A)+Q(B)
-AB$ when $P(x)=Q(x)= gx^4+x^2$. He extended this work
\cite{mehta2,MMN1} with coauthors to matrices coupled in a chain,
model corresponding to $P(A_1,\cdots,A_d)=\sum_{i=1}^d P_i(A_i) 
-\sum_{i=2}^dA_{i-1} A_i$. However, he did not discuss
in these works the limiting spectral distribution of 
the matrices under the corresponding Gibbs measure. On a less rigorous 
ground, P. Zinn Justin  \cite{zinnjustin,zinnjustin1} discussed 
the limiting spectral measures of the matrices following the 
Gibbs measure of the so-called 
Potts model on random graphs, described by $
P(A_1,\cdots,A_d)=\sum_{i=1}^d P_i(A_i) 
-\sum_{i=2}^dA_{1} A_i$. Very interesting work
was also achieved by 
V. Kazakov (in particular for the so-called $ABAB$
interaction case),  A. Migdal and B. Eynard for instance.
We refer to the review \cite{EYNA} of B. Eynard for a general survey.
Matytsin \cite{matytsin}
obtained the first order asymptotics
for spherical integrals, from which he could study the phase
transition of diverse matrix models (see \cite{MZ} for instance).
O. Zeitouni and myself \cite{GZ2} 
gave a complete proof of part of his derivation in \cite{GZ2}
and the present paper is actually finishing to
put his article \cite{matytsin}
on a firm ground.

In this paper, we investigate the problem of
the first order asymptotics 
of matrix integrals with $AB$ interaction, including the above
Ising model, Potts model, matrix model coupled in a chain
and induced QCD models.
The integration will hold over either Hermitian 
matrices or symmetric matrices.
The case of symplectic matrices could be handle similarly.
We obtain, as a consequence of \cite{GZ2},
 the convergence of the free energy
and represent its limit as the solution of a variational 
problem. We here study this variational problem
and characterize its critical points.
One of the main outcome of this study is to
show that under the Gibbs measure $\mu^N_P$
of  the Ising model described by
$$P(A,B)=P(A)+Q(B)-AB$$
with $P(x)\ge ax^4+b$ and $Q(x)\ge cx^4+d$ with $a,b>0$,
the spectral measures of $(A_1^N,A_2^N)$
converges almost surely and to characterize 
its limit.
More precisely,
we shall prove that

\begin{theo}\label{intro}
1) $(\mun_A,\mun_B)$ converges 
almost surely towards a unique couple $(\mu_A,\mu_B)$
of probability measures on $\R$.

2) $(\mu_A,\mu_B)$ are compactly supported
with finite non-commutative entropy $$\Sigma(\mu)=\int\int\log
|x-y|d\mu(x)d\mu(y).$$

3) There exists a couple $(\rho^{A\ra B}, u^{A\ra B})$ 
of measurable functions on $\R\ts (0,1)$ such that
$\rho_t^{A\ra B}(x)dx$ is a probability measure on $\R$
for all $t\in (0,1)$ and
$(\mu_A,\mu_B, \rho^{A\ra B}, u^{A\ra B})$ 
are  characterized uniquely 
as the minimizer of a strictly convex
function under a linear constraint (see Theorem \ref{limitpoints}).

In particular, $(\rho^{A\ra B}, u^{A\ra B})$ 
are solution of the Euler equation
for isentropic flow with negative pressure
$p(\rho)=-{\pi^2\over 3} \rho^3$
such that, 
 for all $(x,t)$ in the
interior of $\Omega=\{(x,t)\in\R\ts [0,1]; \rho^{A\ra B}_t(x)\neq 0\}$,

\begin{equation}\label{euleint}\left\lbc\begin{array}{l}
\partial_t\rho^{A\ra B}_t+\partial_x(\rho^{A\ra B}_t u^{A\ra B}_t)=0\\
\partial_t(\rho^{A\ra B}_t u^{A\ra B}_t) +\partial_x({1\over 2}
\rho^{A\ra B}_t(u^{A\ra B}_t)^2 -{\pi^2\over 3}(\rho^{A\ra B}_t)^3)=0\\
\end{array}\right.
\end{equation}
with the probability 
measure  $\rho^{A\ra B}_t(x)dx$ weakly converging towards $\mu_A(dx)$
(resp. $\mu_B(dx)$) as $t$ goes to zero (resp. one).

\nn
Moreover, we have  
$$P'(x)-x -{\beta\over 2}u^{A\ra B}_0(x)
-{\beta\over 2} H\mu_A(x)=0,\quad
\mbox{ and }
Q'(x)-x+{\beta\over 2}u^{A\ra B}_1(x)-{\beta\over 2} H\mu_B(x)=0
$$
in the sense of distributions that for all $h\in\Ca_b^1(\R)$, 
$${\b\over 4}\int {h(x)-h(y)\over x-y} d\mu_A(x)d\mu_A(y)
= \int (P'(x)-x -{\beta\over 2}u^{A\ra B}_0(x))h(x) d\mu_A(x),$$
$$
{\b\over 4}\int {h(x)-h(y)\over x-y} d\mu_B(x)d\mu_B(y)
= \int (Q'(x)-x +{\beta\over 2}u^{A\ra B}_1(x))h(x) d\mu_B(x).$$
A more detailed characterization of $(\mu_A,\mu_B, \rho^{A\ra B}, 
u^{A\ra B})$  is given in Theorem \ref{limitpoints}.
\end{theo}
Here, $H\mu$ stands for 
the Hilbert transform of the probability measure $\mu$ given by
$$H\mu(x)=PV\int {1\over x-y} d\mu(y)=\lim_{\e\downarrow 0}
\int {(x-y)\over (x-y)^2+\e^2}d\mu(y)$$

To obtain such a result, we shall first study the
limit obtained in \cite{GZ2}
for spherical integrals. This limit was indeed
given by the infimum of a rate function over
measure-valued processes with given initial and terminal data.
We show in section \ref{sph} that this
infimum is in fact taken 
at a unique probability measure-valued path,
 solution of the  Euler equation for isentropic
flow described in (\ref{euleint}). Using a saddle point method,
we derive from \cite{GZ2} in Theorem \ref{lim}  formulae for
the limiting 
 free energy of some  matrix models with $AB$ 
interaction. In the Ising model case, this
free energy is indeed written as the infimum 
of a strictly convex function,
from which uniqueness of the minimizers
is obtained. As a consequence, we obtain
the convergence of the spectral measures 
under the  Gibbs
measure for Ising model. A variational study then shows 
that the limiting spectral measures satisfies the
above set of equations (see Theorem \ref{limitpoints}) .
For the other considered models ($q$-Potts model,
matrix coupled in a chain, induced QCD), obvious convexity
arguments and therefore uniqueness 
is lost in general, but still holds in certain cases. However, we can still
specify some properties of the limit points (see Theorem 
\ref{limitpoints2}).

In this paper, we shall denote $\Ca([0,1],\Pa(\R))$
the set of continuous processes with values in 
the set $\Pa(\R)$ of probability measures on $\R$,
endowed with its usual weak topology.
For a measurable set
$\Omega$ of $\R\ts [0,1]$,
$\Ca^{2,1}_b(\Omega)$  denotes the
set of  real-valued functions on $\Omega$ which are p times continuously 
differentiable with respect to
the (first) space variable 
and $q$ times  continuously 
differentiable with respect to
the (second) time variable  with bounded 
derivatives.
$\Ca^{p,q}_c(\Omega)$ 
will denote the functions
of $\Ca^{p,q}_b(\Omega)$  with compact support 
in the interior of the measurable set $\Omega$. $L_p(d\mu)$ will
denote the space of measurable functions
with finite $p^{th}$ moment under a given
measure $\mu$. We shall say that
an equality holds  in the sense of distribution
on a measurable set $\Omega$ if it holds, once integrated
with respect to
any $\Ca^{\infty,\infty}_c(\Omega)$ functions.
\section{Study of the rate function
governing the asymptotic behaviour
of  spherical integrals}\label{sph}

In \cite{GZ2}, Ofer Zeitouni 
and I  studied the 
so-called spherical integral
$$
I_N^{(\beta)}(D_N,E_N):=\int \exp\{N\tr(UD_NU^*E_N)\} dm_N^\beta(U),$$
where $m_N^\beta$ denotes the Haar measure on the orthogonal group
 $\Oa_N$
when $\beta=1$ and  on the
 unitary group
${\cal U}_N$
when $\beta=2$, and $D_N,E_N$ are diagonal real matrices whose spectral
measures converge to $\mu_D,\mu_E$. We proved (see Theorem 1.1
in \cite{GZ2})
the existence and
represent as solution to a variational problem
the limit 
$$I^{(\beta)}(\mu_D,\mu_E):=\lim_{N\ra\infty}
 N^{-2} \log I_N^{(\beta)}(D_N,E_N).$$
This result 
in fact was obtained under the additionnal technical assumptions
that
there exists 
a compact subset ${\cal K}$ of $\R$  
such that $\mbox{\rm supp}\, \mun_{D_N}
\subset {\cal K}$
for all $N\in\N$ and   that
$\mun_{E_N}(x^2)$ is uniformly bounded (in $N$).
These hypotheses will be made throughout this section.

\medskip

In this section, we investigate the variational problem
which defines $I^{(\beta)}$
and study its minimizer.
We indeed  prove Matytsin's heuristics \cite{matytsin}
outlined in section 6 of \cite{GZ2}.
Let us recall the formula obtained in \cite{GZ2}
for $I^{(\beta)}$ :
$$
I^{(\beta)}(\mu_D,\mu_E):=
-J_\beta(\mu_D,\mu_E)+I_\beta(\mu_E)-\inf_{\mu\in \PP(\R)}
 I_\beta(\mu)+{1\over 2}\int x^2d\mu_D(x)$$
where, for any $\mu\in\PP(\R)$, 
$$I_\beta(\mu)={1\over 2}\int x^2 d\mu(x) -{\beta\over 2}
\Sigma(\mu).
$$
$J_\beta(\mu_D,.)$ is the rate function governing the deviations
of the law of the spectral measure of $X_N=D_N+W_N$
with a Hermitian (resp. symmetric)  Gaussian Wigner matrix $W_N$ 
and 
a deterministic diagonal
matrix $D_N=\mbox{diag}(d_1,\cdots,d_N)$,
$(d_i)_{1\le i\le N}\in\R^N$,
with  spectral measure $\mun_{D_N}=N^{-1}
\sum_{i=1}^N \d_{d_i}$ weakly converging
towards $\mu_D\in\Pa(\R)$. It
is given (see \cite{GZ2}) by
\begin{equation}\label{minus}
J_\beta(\mu_D,\mu)=\frac{\beta}{2}
\inf\{ S_{\mu_D}(\nu_.) ; \nu\in \Ca([0,1],\Pa(\R)) :\nu_1=\mu\}.
\end{equation}
if 
$$S_{\mu_D}(\nu):=\left\lbc
\begin{array}{l}
+\infty\,, \qquad\mbox{ if } \nu_0\neq\mu_D,\\
S^{0,1}(\nu):=\sup_{f\in \CC^{2,1}_b(\R\ts [0,1])}
\sup_{0\le s\le t\le 1}\bar S^{s,t}(\nu,f)\,,
\qquad\mbox{ otherwise. }\\
\end{array}\right.$$
Here,
we have  set, for
any  $f,g\in \CC^{2,1}_b(\R\ts [0,1])$, any $s\le t\in [0,1]$,
and any $\nu_.\in\CC([0,1],\Pa(\R))$,
\begin{eqnarray}
 S^{s,t}(\nu,f)&=&\int f(x,t)d\nu_t(x)-\int
f(x,s)d\nu_s(x)\nonumber\\
&&-\int_s^t\int\partial_uf(x,u)d\nu_u(x)du
-{1\over2}\int_s^t\int\int{\partial_xf(x,u)-\partial_xf(y,u)\over
x-y}d\nu_u(x)d\nu_u(y)du,\label{it1}\\
\nonumber
\end{eqnarray}
\begin{eqnarray}
 <f,g>_{s,t}^\nu &=&\int_s^t\int \partial_xf(x,u)\partial_xg(x,u)
d\nu_u(x)du\,,\label{it2}
\\
\nonumber
\end{eqnarray}
and
\begin{equation}\label{it3}
\bar S^{s,t}(\nu,f)=  S^{s,t}(\nu,f)-{1\over 2} <f,f>^\nu_{s,t}\,.
\end{equation}
It can be shown by Riesz's theorem (see such
a derivation in \cite{CDG1} for instance) that
any measure-valued path $\nu_.\in\CC([0,1],\PP(\R))$
in 
$\{S_{\mu_D}<\infty\}$ is such that there exists
a process $k_.$ so that

\begin{enumerate}
\item 
$$\inf_{f\in \CC^{2,1}_b(\R\ts [0,1])}
<f-k,f-k>_{0,1}^\nu=0$$
\item $\nu_0=\mu_D$ and for any $f\in \CC^{2,1}_b(\R\ts [0,1])$,
any $0\le s\le t\le 1$, 
\begin{equation}\label{eq1}
S^{s,t}(\nu,f)=<f,k>_{s,t}^\nu.\end{equation}
\end{enumerate}
Then, it is not hard to show that 
$$S_{\mu_D}(\nu_.)={1\over 2}<k,k>_{0,1}^\nu.$$
Therefore,
$J_\beta(\mu_D,\mu)$ is given also by
\begin{equation}\label{eq2}
J_\beta(\mu_D,\mu)=\frac{\beta}{4}
\inf\{ <k,k>_{0,1}^\nu ;
\quad (\nu,k) \mbox{ satisfies }(C)
\}.
\end{equation}
with (C) the condition

{\it (C) :
$\nu_0=\mu_D$, $\nu_1=\mu$, $\partial_x k\in
\overline{\Ca^{1,1}(\R\ts[0,1]))}^{L^2(d\nu_t dt)}$, and, for any 
$f\in \CC^{2,1}_b(\R\ts [0,1])$, any $s,t\in [0,1]$,

$$S^{s,t}(\nu,f)=<f,k>_{s,t}^\nu$$
}

\nn
The main Theorem
of this section
states as follows
\begin{theo}\label{mintheo}
Let   $\mu_E\in\{ J_\beta(\mu_D,.)<\infty\}$
with finite entropy $\Sigma$.
 Then,
the infimum 
in $J_\beta(\mu_D,\mu_E)$ is reached at
a unique probability measures-valued path
$\mu^*\in\Ca([0,1],\Pa(\R))$
such that

\begin{itemize}
\item $\mu_0^*=\mu_D,\ \mu^*_1=\mu_E$.
\item For any $t\in (0,1)$,
$\mu^*_t$ is absolutely continuous with respect to
Lebesgue measure ; $\mu^*_t(dx)=\rho^*_t(x)dx$. $t\in [0,1]\ra \mu^*_t\in\Pa(\R)$
is continuous and therefore
$\lim_{t\downarrow 0} \mu^*_t=\mu_D$, $\lim_{t\uparrow 1}
 \mu^*_t=\mu_E$.
\item Let $k^*$ be such that 
the couple $(\mu^*,k^*)$ satisfies
(C). Then, if we set
$$u^*_t=\partial_x k^*_t+H\mu^*_t(y),$$
$(\rho^*,u^*)$ satisfies the
Euler equation
for isentropic flow described by the
equations, 
for $t\in (0,1)$,  

\begin{eqnarray}
\partial_t \rho^*_t(x)&=&-\partial_x( \rho^*_t(x)u^*_t(x))
\label{eulpr10}\\
\partial_t(\rho^*_t(x)u^*_t(x))&=&-\partial_x({1\over 2}\rho^*_t(x)
u^*_t(x)^2-{\pi^2\over 3}\rho^*_t(x)^3)\label{eulpr20}\\
\nonumber
\end{eqnarray}
in the sense of distributions that
for all $f\in \Ca_c^{\infty,\infty}({\R}\ts [0,1])$,
$$\int_0^1\int\partial_t f(t,x)d\mu^*_t(x)
dt +\int_0^1\int\partial_x  f(t,x) u^*_t(x) d\mu^*_t(x)dt=0$$
and, for any $\e>0$,
any $f\in\Ca^{\infty,\infty}_c(\Omega_\e)$ with
$\Omega_\e:=\{(x,t)\in\R\ts[0,1] :
\rho^*_t(x)>\e\}$,
\begin{equation}\label{distf}\int \left(2u^*_t(x)\partial_t f(x,t)+\left( u^*_t(x)^2
-\pi^2 \rho^*_t(x)^2\right) \partial_x f(x,t)\right)dxdt=0.
\end{equation}
\end{itemize}
If we assume that 
$(\mu_D,\mu_E)$ are  compactly
supported
probability measures, we additionnally
know that $(\rho^*,u^*)$ are smooth in the interior of $\Omega_0$,
which guarantees that (\ref{eulpr10}) and (\ref{eulpr20})
hold everywhere in the interior of $\Omega_0$. Moreover,
$\Omega_0$ is bounded in $\R\ts [0,1]$.
Furthermore,
 there exists a sequence $(\phi^\e)_{\e>0}$ of functions 
such that 
$$\partial_t\phi^\e_t(x) +\left({\partial_x\phi^\e_t(x)\over
2}\right)^2\ge 0$$
and if we set 
$$\rho^\e_t(x):=\pi^{-1} 
( \partial_t\phi^\e+4^{-1}(\partial_x\phi^\e)^2 )^{1\over 2}$$
then

$$
\int (u^*_t(x)-\partial_x{\phi^\e_t(x)\over 2}
)^2 d\mu^*_t(x)dt+ {\pi^2\over 3}\int \left( \rho^*_t(x)-
\rho^\e_t(x)\right)^2 
\left(\rho^*_t(x)+\rho^\e_t(x) \right)dxdt\le\e.$$
\end{theo} 
As a consequence, if we 
let   $\Pi^*_t(x)=\int^x u^*_t(y)dy$, 
which should be thought as the limit in $H_1(\rho^*_t(x)dxdt)$
of the sequence
$(2^{-1}\phi^\e)_{\e>0}$,
we find that it satisfies,  in the sense of distributions
in $\Omega_0$,

$$\partial_t\Pi_t^*=-{1\over2}(\partial_x\Pi^*_t)^2 
 +{\pi^2\over 2}(\rho_t^*)^2,$$
which is Matytsin 's equation \cite{matytsin}.

The (non trivial) existence of 
solutions to the Euler equation for isentropic flow (\ref{eulpr10}),
(\ref{eulpr20}),
 is a consequence of
our variational study. The uniqueness of the
solutions to these equations
 could be derived, under some additional 
regularity properties,  from a convexity 
property of our rate function $S_{\mu_D}$. Even
when such solutions are not unique, we know that
our minimizer is unique due  to a convex
property of $S_{\mu_D}$ which  is a consequence of
 its representation of Property \ref{formula}.1)
below (see Property \ref{unique}).

\begin{pr}\label{formula}Let $\mu_E\in\{ J_\beta(\mu_D,.)<\infty\}$
having finite entropy $\Sigma$. Then, 

1) For any $\nu\in \Ca([0,1],\Pa(\R))$, if $(\nu,k)$ verifies
(C) and $u_t(x)=\partial_x k_t(x)+H\nu_t(x)$,
$$S_{\mu_0}(\nu)={1\over 2}\int_0^1\int (u_t(x))^2
d\nu_t(x) dt+{1\over 2}\int_0^1\int (H\nu_t(x))^2
d\nu_t(x) dt-{1\over 2}(\Sigma(\mu_1)-\Sigma(\mu_0)).$$

%
%

2) Consequently, we can  write $J_\beta$ under the
following form 
\begin{eqnarray}
J_\beta(\mu_D,\mu_E)&=&{\beta\over 4}
\left\lbc \int_0^1\int (u^*_t(x))^2 d\mu^*_t(x) dt
+\int_0^1\int (H\mu^*_t(x))^2 d\mu^*_t(x) dt 
-(\Sigma(\mu_E)-\Sigma(\mu_D))\right\rbc\label{fJ}\\
\nonumber
\end{eqnarray}

with $(\mu^*,u^*)$ as in   Theorem
\ref{mintheo}. Note
here that $\mu^*_t(dx)=\rho^*_t(x)dx$ for $t\in (0,1)$ and
 $\rho^*_.\in L_{3}(dx dt)$, so that
$$\int_0^1\int (H\mu^*_t(x))^2 d\mu^*_t 
dt ={\pi^2\over 3}\int_0^1\int (\rho^*_t(x))^3
dx dt.$$

3) As a consequence,

\begin{eqnarray}
I^{(\beta)}(\mu_D,\mu_E)&=&
-{\beta\over 4}
\left\lbc \int_0^1\int (u^*_t(x))^2 d\mu^*_t (x)dt
+\int_0^1\int (H\mu^*_t(x))^2 d\mu^*_t(x) dt \right\rbc\nonumber\\
&&-{\beta\over 4}(\Sigma(\mu_E)+\Sigma(\mu_D))
+{1\over 2} \int x^2d\mu_D(x)+{1\over 2} \int x^2d\mu_E(x)
-\inf I_\beta.
\label{fI}\\
\nonumber
\end{eqnarray}

\end{pr}
In \cite{matytsin}, a similar result
was announced  (see formulae (1.4) and (2.8) of \cite{matytsin}).
However, it seems (as far as I could understand)
 that in  formulae (\ref{fJ},\ref{fI})  of \cite{matytsin},
the first term as the opposite sign. 
But, in \cite{MZ}, formula (2.18),
the very same result is stated.

\medskip


\medskip

Let us also notice 
that the minimizer $\mu^*_t$ 
has the following representation
in the free probability context.
Let $(\Aa,\tau)$ be a non-commutative 
probability space on which an operator
 $D$ with distribution $\mu_D$,
an operator $E$ with distribution $\mu_E$ 
and a semi-circular variable $S$, free
with $(D,E)$, live. Then, there exists a
joint distribution of $(D,E)$ such
that $(\mu^*_t)_{t\in [0,1]}$
is
the law of a free Brownian bridge
$$X_t=tE+(1-t)D+\sqrt{t(1-t)}S.$$
 The isentropic 
Euler equation which
governs $\mu^*$ hence
partially specify 
the joint law of $(D,E)$.
More specifically, for any $t\in (0,1)$,
our result implies that for any $p\in\N$,
$$\lim_{N\ra\infty}
\int \left(tUE_NU^*+(1-t)D_N\right)^p{e^{N\tr(UE_NU^* D_N)}\over
I_N^{(\beta)}(D_N,E_N)}dm_N^\beta(U)=\tau(
(tE+(1-t)D)^p)=\int x^pd\nu^*_t(x)
$$
if 
$\nu^*_t$ is the unique compactly supported 
probability measure such that, if $S$ is a semicircular variable
free with $(D,E)$, for any $p\in\N$,
$$\int x^pd\mu^*_t(x)=\tau \left( (tE+(1-t)D+\sqrt{t(1-t)}S)^p\right)=
\int x^p d\nu^*_t{\tiny\boxplus}\s_{t(1-t)}(x)$$
where $\boxplus$ denotes the free convolution and $\s_\d$ the
semicircular variable with covariance $\d$.

\subsection{Study of $S_{\mu_0}$}
Hereafter and to simplify the notations,
$\mu_D=\mu_0$ and $\mu_E=\mu_1$
with some probability measures $(\mu_0,\mu_1)$ on $\R$. 
We shall in this section study the rate function $S_{\mu_0}$
and show that it achieves its minimal value on
$\{\nu\in\Ca([0,1],\Pa({\R})):\nu_1=\mu_1\}$ at a
unique continuous measure-valued path $\mu^*$. 

\subsubsection{$S_{\mu_0}$ achieves its minimal value}
Recall that for any probability measure $\mu_0\in\Pa({\R})$,
$S_{\mu_0}$ is a good rate function on $\Ca([0,1],\Pa({\R}))$
(see Theorem 2.4(1) of \cite{GZ2}).
Therefore, the infimum defining $J_\beta(\mu_0,\mu_1)$
is, when it is finite, achieved in $\Ca([0,1],\Pa({\R}))$.
We shall in the sequel restrict ourselves
to $(\mu_0,\mu_1)$ such
that $J_\beta(\mu_0,\mu_1)$ is finite.

\smallskip

\subsubsection{A new formula for $S_{\mu_0}$}\label{newf}

In this section, we shall give a simple 
formula of $S_{\mu_0}(\nu)$ in terms of
$u_.=H\nu_.+\partial_x k_.$ and $\nu$
when $(k,\nu)$ satisfies (C).
We begin with the following 
preliminary  Lemma

\begin{lem}\label{l3}
Let $(\mu_0,\mu_1)\in \{\mu\in\Pa({\R}):\Sigma(\mu)>-\infty\}$
and  $\nu_.\in\Ca([0,1],\Pa({\R}))$
such that $\nu_0=\mu_0,\ \nu_1=\mu_1$
and  $\nu_.\in \{S_{\mu_0}<\infty\}$.
Then, for almost all  $t\in (0,1)$, $\nu_t(dx)\ll dx$
and $$\int_0^1\int
\left(H\nu_t(x)\right)^2
d\nu_t(x)dt
={\pi^2\over 3} \int_0^1 \int \left({d\nu_t(x)\over dx}\right)^3 dx dt
<\infty.$$
\end{lem}
The idea of the proof of the
lemma  is quite simple ; we make, in the definition of $S_{\mu_0}$,
the change of variable $f(x,t)\ra  f(x,t)-\int\log|x-y|d\nu_t(x).$
However, because $(x,t)\ra\int\log|x-y|d\nu_t(x)$
is not in $\Ca^{2,1}({\R}\ts [0,1])$ in general,
the full proof requires approximations of the path $\nu_.$ and becomes
rather technical.
This is the reason why I defer it to
the appendix,
section \ref{app2}.
We shall now prove the following
\begin{pr}\label{nfS}
Let $(\mu_0,\mu_1)\in \{\mu\in\Pa({\R}):\Sigma(\mu)>-\infty\}$
and  $\nu_.\in\Ca([0,1],\Pa({\R}))$
such that $\nu_0=\mu_0,\ \nu_1=\mu_1$
and  $\nu_.\in \{S_{\mu_0}<\infty\}$.
Then,
if $(\nu,k)$ satisfies (C) and if we set
$u_t:=\partial_x k_t(x)+H\nu_t(x),$
we have

$$S_{\mu_0}(\nu)={1\over 2}\int_0^1\int (u_t(x))^2
d\nu_t(x) dt+{1\over 2}\int_0^1\int (H\nu_t(x))^2
d\nu_t(x) dt-{1\over 2}(\Sigma(\mu_1)-\Sigma(\mu_0)).$$
\end{pr}
\prf 

Let us recall that $(\nu,k)$ satisfying
condition (C) implies that 
 for any 
$f\in \CC^{2,1}_b({\R}\ts [0,1])$,

\begin{eqnarray}
\int f(x,t)d\nu_t(x)-\int
f(x,s)d\nu_s(x)&=& 
\int_s^t\int\partial_vf(x,v)d\nu_v(x)ds\nonumber\\
&&
+{1\over2}\int_s^t\int\int{\partial_xf(x,v)-\partial_xf(y,v)\over
x-y}d\nu_v(x)d\nu_v(y)dv\nonumber\\
&&
+ \int_s^t\int \partial_xf(x,s)\partial_x 
k(x,s) d\nu_s(x)
\label{it12}\\
\nonumber
\end{eqnarray}
with $\partial_x 
k\in L^2(d\nu_t(x)\ts dt)$.
Observe that by \cite{tricomi}, p. 170,
for any $s\in [0,1]$ such that
 $\nu_s$ is absolutely continuous 
with respect to
Lebesgue measure with density $\rho_s
\in L_3(dx)$, for any compactly supported 
measurable function $\partial_x f(.,s)$,
$$\int\int{\partial_xf(x,s)-\partial_xf(y,s)\over
x-y}d\nu_s(x)d\nu_s(y)=2\int \partial_xf(x,s) H\nu_s(x) dx ds.$$
Since by  Lemma
\ref{l3}, for almost all $s\in [0,1]$,  $\nu_s(dx)\ll dx$
with a density $\rho_s\in L_3(dx)$
 we conclude that, 
 in the sense of distributions on $\R\ts [0,1]$, (\ref{it12}) implies
\begin{equation}\label{mn}
\partial_s\rho_s+\partial_x(u_s \rho_s)=0,
\end{equation}
i.e for any compactly supported
$f\in\Ca^{\infty,\infty}_c({\R}\ts [0,1])$ 
vanishing at the boundary of ${\R}\ts [0,1]$, 
$$\int_0^1\int_{\R} (\partial_s f(x,s)+ u_s
\partial_x f(x,s))\rho_s(x)dx =0.$$
Note here that, by dominated convergence theorem, 
we can equivalently take $f\in\Ca^{2,1}_b({\R}\ts [0,1])$.

Moreover, since 
$H\nu_.$ belongs to $L^2(d\nu_s\ts ds)$ by  Lemma \ref{l3},
we can write
\begin{eqnarray}
2S_{\mu_0}(\nu_.)=<k,k>_{0,1}^\nu &=& 
\int_0^1\int_{{\R}} (u_s(x))^2 d\nu_s(x)ds
+\int_0^1\int_{{\R}} (H\nu_s(x))^2 d\nu_s(x)ds\nonumber\\
&&
-2\int_0^1 \int_{{\R} }
H\nu_s(x)u_s(x) d\nu_s(x)ds\label{as}\\
\nonumber
\end{eqnarray}
We shall now see that the last 
term in the above right hand side only
depends on $(\mu_0,\mu_1)$. The only difficulty 
in the proof of this point lies in the fact
that $x,s\in {\R}\ts[0,1] \ra \int\log|x-y|d\nu_s(y)$
is not in $ \Ca^{2,1}_c({\R}\ts [0,1])$.

However, following  Lemma 5.16 in \cite{CDG2},
 if ${\tiny\boxplus}$ denotes 
the free convolution (see \cite{Vstflour} for a definition),
if  for any $\d>0$, 
$\s_\d$ denotes
the semicircular law with covariance $\d$,
and if $u^\d_t$ denotes the field 
corresponding to $\nu_t{\tiny\boxplus}\s_\d$,

\begin{equation}\label{asd}
\Sigma(\nu_1{\tiny{\tiny\boxplus}}\s_\d)-\Sigma(\nu_0{\tiny\boxplus}\s_\d)
=2\int_0^1 \int_{{\R} }
H(\nu_s{\tiny\boxplus}\s_\d)(x) u_s^\d(x) d\nu_s
{\tiny\boxplus}\s_\d(x)ds.
\end{equation}
It is well known that $\d\ra \Sigma(\nu{\tiny\boxplus}\s_\d)$
is continuous (see \cite{SV}, Theorem 2.1
for the lower semicontinuity and use the
well known upper semi-continuity). 
Moreover, if $X_s$ is a random variable with
distribution $\nu_s$ and $S$ a semicircular variable,
free with $X_s$,  living in a non commutative
probability space $(\Aa,\tau)$,
by Theorem 4.2 in \cite{CDG2}, the field 
$u^\d$  is given,
$\nu_s{\tiny\boxplus}\s_\d$ almost surely, by 
$$ u^\d_s
=\tau(\partial_xk(X_s ,s)|X_s+\sqrt{\d}S)+
H\nu_s{\tiny\boxplus}\s_\d.
$$ 
Consequently,
$$\int_{{\R}} 
H(\nu_s{\tiny\boxplus}\s_\d)(x)u^\d_s(x) d\nu_s
{\tiny\boxplus}\s_\d(x)=\tau\left( \partial_xk(X_s,s)
H(\nu_s{\tiny\boxplus}\s_\d)(X_s+\sqrt{\d}S)\right)
+\tau\left(H(\nu_s{\tiny\boxplus}\s_\d)(X_s+\sqrt{\d}S)^2\right) .$$
Moreover, by Voiculescu \cite{Vo5}, Proposition 3.5
and
Corollary 6.13,
if  $\nu_s(dx)=\rho_s(x)dx\in
L_3(dx)$,

$$\lim_{\d\ra 0} \tau\left( (H(\nu_s{\tiny\boxplus}\s_\d)(X_s+\sqrt{\d}S)
-H\nu_s(X_s))^2\right)=0.$$
Therefore, for any such $s\in [0,1]$,
\begin{equation}\label{cvb}
\lim_{\d\ra 0}\int_{{\R}} 
H(\nu_s{\tiny\boxplus}\s_\d)(x)u^\d_s(x)
d\nu_s
{\tiny\boxplus}\s_\d(x)
= \int_{{\R} }
H\nu_s(x)u_s(x) d\nu_s(x).
\end{equation}
Note that 
by  Lemma \ref{l3}, this convergence holds for
almost all $s\in [0,1]$ since $\rho_.\in L_3(dxdt)$.
Finally, by Propositions 3.5 and  3.7 of \cite{Vo5},
for any $s$ such that $H\nu_s$ is well defined,

$$H \nu_s{\tiny\boxplus}\s_\d(X_s+\sqrt{\d}S)=
\tau( H\nu_s(X_s)|X_s+\sqrt{\d}S)$$
so that for any $\d>0$
$$\tau \left((H \nu_s{\tiny\boxplus}\s_\d(X_s+\sqrt{\d}S))^2\right)
\le \tau\left((H \nu_s(X_s))^2\right).$$
Therefore, dominated convergence theorem and 
(\ref{cvb}) imply that

$$\lim_{\d\ra 0}\int_0^1 
\int_{{\R}} 
H\nu_s{\tiny\boxplus}\s_\d(x)u_s^\d(x) d\nu_s
{\tiny\boxplus}\s_\d(x) ds=\int_0^1 
\int_{{\R}} 
H\nu_s(x)u_s(x) d\nu_s(x) ds.$$
Thus, (\ref{asd}) extends to $\d=0$
which proves, with (\ref{as}), Property \ref{nfS}.

\hfill\xx

\subsubsection{Uniqueness of the minimizers of $S_{\mu_0}$}

We shall use the 
formula for $S_{\mu_0}$
obtained in the last section
to prove that

\begin{pr}\label{unique}
For any $(\mu_0,\mu_1)\in
\Pa({\R})$ with finite entropy $\Sigma$,
there exists a 
 unique measures-valued path $\mu^*$ such
that
$$J_\beta(\mu_0,\mu_1)={\beta\over 2}\inf\{ S_{\mu_0}(\nu_.)\quad
:\quad
\nu_1=\mu_1\}={\beta\over 2}S_{\mu_0}(\mu^*).$$
\end{pr}
In the following, $\mu^*$ shall always denote 
the minimizer of Property \ref{unique}
and $\partial_x k^*, u^*$ its associated fields.

\nn
\prf
According to the previous section,
the minimizers of  $S_{\mu_0}$ also
minimize

$$S(u,\rho)= {\pi^2\over 3}\int_0^1\int_{{\R}} (\rho_t(x))^3 dxdt
+\int_0^1\int_{{\R}} (u_t(x))^2 \rho_t(x) dxdt$$
under the constraint $\partial_t\rho_t+\partial_x(\rho_t u_t)=0$
in the sense of distributions, $\rho_t\ge 0$
almost surely w.r.t Lebesgue measure and $\int \rho_t(x)dx=1$,
and with given initial and terminal data for $\rho$  given by
$$\lim_{t\downarrow 0}
 \rho_t(x) dx = \mu_0(dx),\quad \lim_{t\uparrow 1}
 \rho_t(x) dx=\mu_1(dx)$$
where convergence holds in the weak sense (with
respect to bounded continuous functions)
and is simply due to
the fact that $S_{\mu_0}$ is finite only on 
 continuous measure-valued paths.

Let   $m=u\rho$ be the corresponding 
momentum.
In the variables $(m,\rho)$, $S(\rho,m)$
reads
$$\Sa(m,\rho)=
{\pi^2\over 3}\int_0^1\int_{{\R}} (\rho_t(x))^3 dxdt
+\int_0^1\int_{{\R}} {(m_t(x))^2 \over \rho_t(x)} dxdt$$
with the convention ${0\over 0}=0$,
whereas the constraint becomes {\bf linear }

$$\partial_t(\rho_t(x))+\partial_x(m_t(x))=0,\quad
\rho_t(x)dx\in\Pa(\R)\quad\forall t\in [0,1], 
\quad \lim_{t\downarrow 0}
 \rho_t(x) dx =\mu_0(dx),\quad \lim_{t\uparrow 1}
 \rho_t(x) dx=\mu_1(dx).
$$
We now observe that $\Sa$
is a {\bf strictly  convex } function. Indeed,
if $(m^1,\rho^1)$ and $(m^2,\rho^2)$ are any two couples
of measurable functions in $\{\Sa<\infty\}$,
it is easy to see that for any $\a\in (0,1)$
\begin{eqnarray*}
\partial_\alpha^2 \Sa( \a m^1+(1-\a)m^2,
\a \rho^1+(1-\a)\rho^2)&=&
2\pi^2\int_0^1\int_{{\R}} (\rho_t^1(x)-\rho_t^2(x))^2
(\a \rho^1_t(x)+(1-\a)\rho^2_t(x))
 dxdt\\
&&
+2\int_0^1\int_{{\R}} {(\rho_t^1(x)m^2_t(x)-\rho_t^2(x)m^1_t(x))^2
\over (\a \rho^1_t(x)+(1-\a)\rho^2_t(x))^3}
 dxdt.\\
\end{eqnarray*}
Hence, 
$\partial_\alpha^2 \Sa( \a m^1+(1-\a)m^2,
\a \rho^1+(1-\a)\rho^2)>0$ for some
 $\a\in (0,1)$ unless  for almost all
$t\in [0,1]$
$$\rho_t^1(x)=\rho_t^2(x) =\rho_t(x),\mbox{ and }
u^1_t(x)={m^1_t(x)\over \rho_t^1(x)}=
{m^2_t(x)\over \rho_t^2(x)}=u^2_t(x)\quad \rho_t(x)dxdt\mbox{ a.s.}$$
In other words,  $\Sa$ is strictly convex. 
By standard convex analysis, the strict convexity of $\Sa$
results with the uniqueness of its minimizers
given a linear constraint, and in particular
in $J_\beta$.  More precisely, from the above, the minimizer $\mu^*$
in $J_\beta$ is defined uniquely for almost all $t\in [0,1]$
(and then everywhere by continuity of $\mu^*$)
and its field $u^*$, or
equivalently $\partial_x k^*$, is then defined  uniquely
$\mu^*_t(dx) dt$
almost surely.

\hfill\xx

\subsection{ A priori properties of the minimizer $\mu^*$ }\label{reg}

In this section, we shall see that the minimizer $\mu^*$
has to be the distribution of a free Brownian bridge
when at least one of the probability measure $\mu_0$ or $\mu_1$
are compactly supported, the other having finite variance
(since we rely on \cite{GZ2}'s results).
To simplify the statements, we shall assume throughout this section
that both probability measures are compactly supported.
This property will unable us to 
obtain a priori properties on the 
laws of the minimizers, such as existence, boundedness,
 and smoothness of
their densities. 

\subsubsection{Free Brownian bridge characterization of the minimizer}

Let us state more precisely the 
theorem obtained in this section. 
A free Brownian bridge
between $\mu_0$ and $\mu_1$ 
is the law of
\begin{equation}\label{fb}
X_t=(1-t)X_0+tX_1+\sqrt{t(1-t)} S
\end{equation}
with a semicircular variable $S$, free with $X_0$ and $X_1$, 
with law $\mu_0$ and $\mu_1$ respectively.
We let ${\mbox{FBB}(\mu_0,\mu_1)}\subset\Ca([0,1],\Pa({\R}))$ 
denote the set of such laws (which depend
of course not only on $\mu_0,\mu_1$ but on
the joint distribution of $(X_0,X_1)$
too). 
Then, we shall prove
that
\begin{theo}\label{fbb} Assume $\mu_0,\mu_1$ compactly supported. 
Then,
\begin{eqnarray*}
J_\beta(\mu_0,\mu_1)&=& {\beta\over
2} \inf\{ S(\nu),\nu_0=\mu_0,\nu_1=\mu_1\}\\
&=& {\beta\over
2} \inf\{ S(\nu)\quad ;\quad \nu\in \mbox{FBB}(\mu_0,\mu_1)
\}.\\
\end{eqnarray*}
Therefore, since $\mbox{FBB}(\mu_0,\mu_1)$ is a closed
subset of $\Ca([0,1],\Pa(\R))$,  the unique minimizer $\mu^*$
in the above infimum belongs to 
 ${\mbox{FBB}(\mu_0,\mu_1)}$.
\end{theo}
The proof of Theorem \ref{fbb}
is rather technical and goes back through the
large random matrices origin of $J_\beta$.
We therefore defer it to the appendix.

\subsubsection{Properties of the free Brownian motion paths}\label{bbp}
As a consequence of  Theorem \ref{fbb},
we shall prove that

\begin{cor}\label{pr}
Assume $\mu_0$ and $\mu_1$ 
compactly supported. Then, 

a) There exists a compact set $K\subset{\R}$
so that for all $t\in [0,1]$, $\mu_t^*(K^c)=0$. For all $t\in (0,1)$,
the support of $\mu^*_t$ is the closure of its interior.

b)
$\mu^*_t(dx)\ll dx$ for all $t\in [0,1]$.
Let $\rho^*_t(x)={d \mu^*_t(x)\over dx}$.

c)
There exists a finite constant $C$ (independent of $t$)
so that, 
$\mu^*_t$ almost surely,
$$\rho^*_t(x)^2 +(H\mu^*_t(x))^2\le (t(1-t))^{-1}$$

and
$$|u^*_t(x)|\le C(t(1-t))^{-{1\over 2}}.$$

d)
$(\rho^*,u^*)$
are analytic in  the interior 
of $\Omega=\{x,t\in\R\ts[0,1] : \rho^*_t(x)>0\}$.

e)
At the boundary of $\Omega_t=\{x\in\R : \rho^*_t(x)>0\}$, 
for $x\in\Omega_t
$,
$$|\rho_t^*(x)^2\partial_x
\rho_t^*(x)|\le {1\over 4\pi^3 t^2(1-t)^2}\quad\Rightarrow
\quad \rho_t^*(x)\le \left({3\over 4\pi^3 t^2(1-t)^2}\right)^{1\over
3}
(x-x_0)^{1\over 3}$$

if $x_0$ is the nearest point of $x$ in $\Omega_t^c$.
\end{cor}
Consequently, the minimizer $\mu^*$ 
may only have shocks at the boundary 
of its support.

\nn
\prf This corollary
is a direct consequence 
of Theorem \ref{fbb} and we shall 
collect these properties for any free brownian bridge 
law. Indeed,
let $(\Aa,\tau)$ be a non-commutative probability 
space in which two operators $X_0,X_1$ 
with laws $\mu_0$ and $\mu_1$ and
 a semicircular variable $S$, free with $(X_0,X_1)$,
live. We assume throughout that
  $X_0$ and $X_1$ are bounded by $C$ for the operator norm (i.e
$\mu_0([-C,C]^c)=\mu_1([-C,C]^c)=0$).

Let $\mu_t$ be  the distribution of 
$$
X_t=tX_1+(1-t)X_0+\sqrt{t(1-t)}S.$$
Clearly, since $S$ is bounded by $2$ for the operator norm, 
$X_t$ is bounded 
by $C+2$ for all $t\in [0,1]$. Thus, proposition 
4 in \cite{Biane}  finishes the proof of  a). 
Following Voiculescu (see  Proposition 3.5
and Corollary  3.9 in \cite{Vo5}),
the Hilbert transform of $\mu_t$ is given,
$\mu_t$-almost surely,  by
$$H\mu_t(x)= \tau( { (2\sqrt{t(1-t)})}^{-1}S
|X_t) $$
with $\tau(\quad|X_t)$ the
conditionnal expectation with
respect to $X_t$, i.e the
orthogonal projection on 
the sigma algebra generated by $X_t$. We deduce that since $S$ is 
bounded for the operator norm by $2$,  $\mu_t$-almost surely,

$$|H\mu_t(x)|\le {1\over \sqrt{t(1-t)}}.$$
Further, following \cite{BCG},
the stochastic differential equation satisfied 
by $X_t$ shows
that, for any twice continuously differentiable 
function $f$ on ${\R}$,
\begin{equation}\label{pol}
\mu_t(f)=\mu_0(f)+{1\over 2}\int_0^t\int\int
{\partial_x f(x)-\partial_x f(y)\over x-y}d\mu_s(x)d\mu_s(y) ds
+\int_0^t\int \partial_x f(x) \partial_x k_s(x)d\mu_s(x)ds
\end{equation}
with $k$ the element of $L^2(d\mu_s ds)$
given by

$$\partial_x k_s(x)=\tau( {X_s-X_1\over s-1}|X_s).$$
Hence,
\begin{equation}\label{baba}
u_t:=\partial_x k_t+H\mu_t= 
\tau( {X_t-X_1\over t-1}|X_t)
+\tau( \sqrt{ 4t(1-t)}^{-1}S
|X_t)= \tau(X_1-X_0+  {(1-2t)\over  2\sqrt{ t(1-t)}} S|X_t).
\end{equation}
Therefore, 
 $\mu_t$-almost surely,
\begin{equation}\label{b1}
|u_t|\le 2C +{ 1\over \sqrt{t(1-t)}
}
.\end{equation}
Moreover,  by Biane's  results \cite{Biane}, we know
that, for $t\in (0,1)$,
 $\mu_t$ is absolutely continuous with respect to Lebesgue
measure. We denote by $\rho_t$ its density. Then,
we also know that for all $t\in (0,1)$,
$\mu_t$-almost surely, 
\begin{equation}\label{b2}
\rho_t(x)^2+(H\mu_t)^2(x)\le {1\over t(1-t)}.\end{equation}
Let us 
mention the  regularity properties that
$(\mu_t)_{t\in (0,1)}$ will inherite from
its free Brownian bridge  formula. 
If $\nu_t$ denotes the law of 
$tX_1+(1-t)X_0$,
we have, following Biane \cite{Biane}, corollary 3, that
if we set
\begin{eqnarray*}
v(u,t)&=&\inf\{ v\ge 0 |
\int \frac{d\nu_t(x)}{(u-x)^2+v^2}\le ({t(1-t)})^{-1}\},
\\
&=& \inf\{ v\ge 0 |
\tau \left( ((tX_1+(1-t)X_0-u)^2+v^2)^{-1} \right) \le
({t(1-t)})^{-1}\},\\
\end{eqnarray*}
$$\psi(u,t)=u+t(1-t)\int \frac{(u-x) d\nu_t(x)}{(u-x)^2+v(u,t)^2},$$
then
$$H{\mu_t}(\psi(u,t))
=\int \frac{(u-x) d\nu_t(x)}{(u-x)^2+v(u,t)^2},$$
while 
$$\rho_t(\psi(u,t))=\frac{ v(u,t)}{\pi t(1-t)}.$$
From these formulae, we observe 
that $\psi^{-1}$ is analytic in the
interior of   $\Omega $ since $\psi'$ 
is bounded below by a 
positive constant there (see Biane \cite{Biane}, p 713
and the obvious analyticity  in the time parameter
$t\in (0,1)$ ), it is clear that $\rho$ 
is $\Ca^\infty$ in $\Omega$. Hence,
the weak equation (\ref{pol}) is verified in the strong
sense in $\Omega$ and 
we find that in $\Omega$,
$u_t(x)=\rho_t(x)^{-1}{\int_x^\infty \partial_t\rho_t(y) dy
}$
is  $\Ca^\infty$. 

At the boundary of $\Omega_t=\{x: (x,t)\in\Omega\}$, Biane (\cite{Biane}, corollary 5)
also noticed 
that
$$|\rho_t(x)^2\partial_x\partial_t(x)|\le
{1\over 4\pi^3 t^2(1-t)^2}\Rightarrow \rho_t(x)\le \left({3\over 4\pi^3
t^2(1-t)^2}\right)^{1\over 3}
(x-x_0)^{1\over 3}$$
with $x_0$ the nearest point of the boundary
of $\Omega_t$ from $x$.

\subsection{The variational problem}

We now turn to
the analysis of the variational
problem defining $J_\beta$ ; we shall prove that

\begin{pr}\label{variational}
Assume that $\mu_0$ and $\mu_1$ are  probability measures 
on $\R$ such that $\Sigma(\mu_0)$ and
$\Sigma(\mu_1)$ are finite.
Then, the path $\mu^*\in \Ca([0,1],\Pa({\R}))$ minimizing 
$J_\beta(\mu_0,\mu_1)
$
satisfies ;

1)$\mu^*_0=\mu_0$ and
$\mu^*_1=\mu_1$. 

2) For any $t\in (0,1)$, $\mu^*_t(dx)\ll dx$. Let $(\rho^*_t)_{t\in (0,1)}$
denote the corresponding density. By continuity of $\mu^*$,
$\mu^*_t(dx)=\rho^*_t(x)dx$ converges towards 
$\mu_0$ (resp. $\mu_1$) as $t$ goes to zero (resp. one)
in the usual weak sense on $\Pa({\R})$.

3)  $\mu^*$ is characterized
as the unique continuous measure-valued
path such that 
$\mu^*_0=\mu_0$ and $\mu^*_1=\mu_1$ 
and, for any $\nu\in \{S_{\mu_0}<\infty\}$ so that $(\nu,k)$
satisfies (C) and $\nu_1=\mu_1$, we have, with $u=\partial_x k+H\nu$,
\begin{eqnarray}
\int [\int 2u^*_t(u_td\nu_t -u_t^*d\mu_t^*)-\int(u^*_t)^2
(d\nu_t-d\mu_t^*)+
\pi^2 \int (\rho^*_t)^2(d\nu_t-d\mu_t^*)]dt\ge 0\label{ineq1p}\\
\nonumber
\end{eqnarray}

4)As a consequence, $(\rho^*,u^*)$ satisfies the
Euler equation
for isentropic flow described by the
equations, 
for $t\in (0,1)$,  

\begin{eqnarray}
\partial_t \rho^*_t(x)&=&-\partial_x( \rho^*_t(x)u^*_t(x))
\label{eulpr1}\\
\partial_t(\rho^*_t(x)u^*_t(x))&=&-\partial_x({1\over 2}\rho^*_t(x)
u^*_t(x)^2-{\pi^2\over 3}\rho^*_t(x)^3)\label{eulpr2}\\
\nonumber
\end{eqnarray}
in the sense of distributions that
for all $f\in \Ca_c^{\infty,\infty}({\R}\ts [0,1])$,
$$\int_0^1\int\partial_t f(t,x)d\mu^*_t(x)
dt +\int_0^1\int\partial_x  f(t,x) u^*_t(x) d\mu^*_t(x)dt=0$$
and, for any $\e>0$,
any $f\in\Ca^{\infty,\infty}_c(\Omega_\e)$ with
$\Omega_\e:=\{(x,t)\in\R\ts[0,1] :
\rho^*_t(x)>\e\}$,
\begin{equation}\label{eulpr3}\int \left(2u^*_t(x)\partial_t f(x,t)+\left( u^*_t(x)^2
-\pi^2 \rho^*_t(x)^2\right) \partial_x f(x,t)\right)dxdt=0.
\end{equation}

\medskip

\nn
Let us now assume that $(\mu_0,\mu_1)$ are compactly supported.
Then

5)  (\ref{eulpr2}) is true everywhere in the interior of
$\Omega_0$. Moreover,  (\ref{eulpr3}) can be improved by the statement 
that
$$ \int_0^1\int_{\R} u^*_t(x)
(\partial_t f(t,x) +u^*_t(x)\partial_x f(t,x)) d\mu^*_t(x)dt
={\pi^2\over 3} \int_0^1\int_{\R} (\rho^*_t(x))^3 \partial_x 
f(t,x) dx dt$$
for all $f\in\Ca^{1,1}_b(\R\ts (0,1))$.

6) There exists a sequence $(\phi^\e)_{\e>0}$ of functions 
such that 
$$\partial_t\phi^\e_t(x) +\left({\partial_x\phi^\e_t(x)\over
2}\right)^2\ge 0$$
and if we set 
$$\rho^\e_t(x):=\pi^{-1} 
( \partial_t\phi^\e+4^{-1}(\partial_x\phi^\e)^2 )^{1\over 2}$$
then

$$
\int (u^*_t(x)-\partial_x{\phi^\e_t(x)\over 2}
)^2 d\mu^*_t(x)dt+ {\pi^2\over 3}\int \left( \rho^*_t(x)-
\rho^\e_t(x)\right)^2 
\left(\rho^*_t(x)+\rho^\e_t(x) \right)dxdt\le\e.$$

\end{pr}

\begin{disc}
 Matytsin \cite{matytsin} noticed
that if we set
$$f(x,t)=u_t^*(x)+i\pi\rho_t^*(x),$$
then the Euler equation for isentropic 
flow implies that $f$ the Burgers equation.
Hence, if one assumes that $f$ can be smoothly 
extended to the complex plan, we find
by usual characteristic methods that for $z\in\C$
$$f( f(z,0)t+z,t)=f(z,0)$$
and therefore, setting $G_+(z)=z+f(z,0)$ and $G_-(z)=
z-f(z,1)$, we see that our problem boils down to solve
$$G_+\circ G_-(z)=G_-\circ G_+(z)=z$$
with $\Im(G_+)(x)=\pi \rho_0(x)$ and $\Im(G_-)(x)=-\pi \rho_1(x)$
if $\rho_0$ and $\rho_1$ are the densities of
$\mu_0,\mu_1$ respectively.
This kind of characterization is 
in fact reminiscent to
the description of minimizers 
provided by P. Zinn Justin \cite{zinnjustin1}.
However, such a result would require 
more smoothness of $(\rho^*,u^*)$ than what we 
proved here.

\end{disc}

\nn
{\bf Proof of Property \ref{variational} :}
By property \ref{nfS}, we want to minimize
$$S(\rho, u):=\int_0^1\int (u_t(x))^2 \rho_t(x)dx dt
+{\pi^2\over 3} \int_0^1\int (\rho_t(x))^3 dx dt$$
under 
 the constraint (C') :
$$\partial_t\rho_t+\partial_x(u_t
\rho_t)=0,\quad\lim_{t\downarrow 0}\rho_t(x)dx=\mu_0,\quad
\lim_{t\uparrow 1}\rho_t(x)dx=\mu_1$$
and when $\rho_t(x)dx\in\Pa(\R)$ for all $t\in [0,1]$.
To study the variational problem associated with this energy, 
I know essentially three ways. The first is to make a perturbation 
with respect to the source. This strategy was followed 
by D. Serre in \cite{SE} but applies only when we know
a priori that $(\rho^*,u^*\rho^*)$ are
uniformly bounded. Since this case corresponds 
to the case where $\mu_0,\mu_1$ are compactly supported,
we shall consider it in the second part of the proof.
One can also use   a target type perturbation,
which is a standard perturbation on the space of probability measure,
viewed as a subspace of
the vector space of measures. This method gives (3) in 
Property \ref{variational} as we shall see.
The last way is
to use convex analysis, following for instance Y. Brenier (see 
\cite{brenier}, section 3.2). We shall also detail 
these arguments, since it provides  some approximation property of the field 
$u^*$, as described in Property \ref{variational}.6).

We begin with the target type perturbation.
In the following, we denote $(\rho^*,u^*)$
the minimizer of $S$ under the constraint (C').
Let $(\rho,u)\in \{S<\infty\}$ satisfying the constraint (C').
Then, for any $\a\in [0,1]$,  we set, with $m=\rho u$
and $m^*=\rho^* u^*$,
$$
\rho^\a=(1-\a)\rho^*+\a\rho,\quad  m^\a:=(1-\a)(\rho^*u^*)
 +\a (\rho u):= \rho^\a u^\a,\quad u^\a=(m^\a/\rho^\a). $$
It is then not hard to check that
$S(\rho^\a, u^\a)<\infty$ for all $\a\in [0,1]$. Moreover, 
by the convexity of $\phi:\a\ra (\rho^\a_t(x))^{-1}(m^\a_t(x))^2+3^{-1} \pi^2
(\rho^\a_t(x))^3$ for all admissible $(\rho,m),(\rho',m')$,
we see that $\a^{-1}(\phi(\a)-\phi(0))$ decreases as $\a\ra 0$
showing, by monotone convergence theorem the 
existence of $\partial_\a S(\rho^\a, u^\a)(0^+)$ and 

\begin{eqnarray*}
\partial_\a S(\rho^\a, u^\a)(0^+)&=& 
\int [-(u^*)^2\rho^*- (u^*)^2\rho+2 mu^* +\pi^2 (\rho^*)^2(\rho-\rho^*)]dxdt\\
&=& \int [2u^*(m-m^*)-(u^*)^2 (\rho-\rho^*)+
\pi^2 (\rho^*)^2(\rho-\rho^*)]dxdt.\\
\end{eqnarray*}
Hence,  for any $(\rho,u)\in \{S<\infty\}$, we have
\begin{eqnarray}
\partial_\a S(\rho^\a, u^\a)(0^+)&=& \int [2u^*(m-m^*)-
(u^*)^2 (\rho-\rho^*)+
\pi^2 (\rho^*)^2(\rho-\rho^*)]dxdt\ge 0\label{ineq1}\\
\nonumber
\end{eqnarray}
Reciprocally, since 
$S$ is convex in $(\rho,m)$, we know that
$$S(\rho^\a, u^\a)\ge S(\rho^*, u^*) +\partial_\a S(\rho^\a, u^\a)(0^+)\a
$$ 
so that (\ref{ineq1}) implies that
$S(\rho^\a, u^\a)\ge S(\rho^*, u^*)$ for all $\a\in [0,1]$
and $(\rho,u)\in\{S<\infty\}$. Hence, (\ref{ineq1}) 
characterizes our unique minimizer, which proves Property
\ref{variational}.3). We can apply  this result 
with 
$$\rho=\rho^*+\e\partial_x \phi,\quad  m=m^*-\e\partial_t\phi$$
for some 
$\phi\in\Ca^{1,1}_c(\Omega_\e)$, $\e>0$,
such that $\partial_x\phi(.,0)=\partial_x\phi(.,1)=0$,
insuring that $S(\rho,u)$ has finite entropy.
This yields the second point of Property \ref{variational}.3).
Conditions at the boundary of the support 
can also be deduced from (\ref{ineq1}),
but they are hardly understandable, since the conditions over the
potentials $\phi$ become more stringent.

\bigskip

To prove  the last points of our property
which concerns the case where $(\mu_0,\mu_1)$ are compactly supported,
 we
follow D. Serre \cite{SE} and Y. Brenier \cite{brenier}.

The idea developped in \cite{SE} is basically to
set $a_t(x)=a(t,x)=(\rho_t^*(x), \rho_t^*(x) u_t^*(x))$
so that $\mbox{div}(a_t(x))=0$
and
perturbe $a$ by considering a family 
$$a^g=J_g(a.\nabla_{x,t} h)\circ g=J_g(\rho^*
(\partial_t h+u^* \partial_x h))\circ g $$
with a $\Ca^\infty$ diffeomorphism $g$ of
$Q=[0,1]\ts\R$
with inverse $h=g^{-1}$ and Jacobian $J_g$.
Such an approach yields the Euler's equation 
(\ref{eulpr3}) of Property \ref{variational}
(use the boundedness of $(\rho^*,u^*)$
obtained in Corollary \ref{pr} to apply theorem 2.2 of \cite{SE}).
Moreover, since we saw in Corollary \ref{pr}.d) that
$\rho^*$ and $u^*$ are smooth in the interior
of $\Omega_0$, 
(\ref{eulpr2}) results with Property \ref{variational}.5).

We now 
 developp convex analysis
for our problem following \cite{brenier}.    By Corollary \ref{pr}.a),
we see that there exists a compact $K$ such that $\mu_t^*(K^c)=0$
for all $t\in [0,1]$.
We set $Q=K\ts[0,1]$
and $E=\Ca_b(Q)\ts \Ca_b(Q)$.

For 
 any  continuous 
functions $F,\Phi\in E$ , we set
$$\alpha(F, \Phi)={2\over 3\pi}\int_Q  |F(x,t)+({\Phi(x,t)\over 2})^2|^{3\over
2} dxdt
$$
 if  $ F+({\Phi\over 2})^2\ge 0,$
on $Q$, and $+\infty$ otherwise.  
For any $(\mu,M)\in E'$, let us consider
$$\alpha^*(\mu,M)=
\sup\{ \int_{Q}F(x,t)\mu(dx,dt) +\int_Q 
\Phi(x,t) M(dx,dt)-\alpha(F ,\Phi)\}.$$
It is not hard to see that $\alpha^*(\mu,M)<\infty$ 
iff 
 $\mu$ is  non negative, $M$ is absolutely continuous w.r.t
$\mu$ (to prove  both cases take $F+4^{-1}\Phi^2=0$
and argue by choosing $\Phi$ wisely)
and $\mu$ is absolutely continuous w.r.t
Lebesgue measure with density in $L_3(dxdt)$
(in this case, take $\Phi=0$). Moreover,
if we denote $\mu(dx,dt)=\rho_t(x)dxdt$, $M(dx,dt)=
m_t(x) dxdt$, it is not hard to
see that $\alpha^*(\mu,M)=\Sa(\rho,m)$.
Now, let  
$$\beta(F,\Phi)= \int_Q F(x,t)\rho^*_t(x)dxdt+\int_Q \Phi(x,t) 
u^*_t(x) \rho^*_t(x)dxdt 
$$
if there exists $\phi\in\Ca_b^{1,1}(Q)$
such that
$$F(x,t)+\partial_t\phi(x,t)=0,\quad \Phi(x,t)+\partial_x\phi(x,t)=0$$
for all $(x,t)\in Q$, and is equal to $+\infty$
otherwise. We consider
\begin{eqnarray*}
\beta^*(\mu,M)&=&\sup\{ \int_{Q}F(x,t)\mu(dx,dt) +
\int_Q \Phi(x,t) M(dx,dt)-\beta(F,  \Phi)\}\\
\end{eqnarray*}
Then, $\beta^*$ is infinite
unless 
$\int_Q F(x,t) (\mu(x,t)-\rho^*_t(x))dxdt + \int_Q\Phi(x,t) (M(x,t)
-m^*_t(x))dx dt=0$
for all $(F,\Phi)\in E$. Therefore, $\partial_t \mu+\partial_xM=0$
in the sense of distributions, 
 $\int \mu(x,t)dx=1$ for almost all
$t\in [0,1]$ and $\lim_{t\downarrow 0}
\mu(x,t)dx =d\mu_0(x)$, $\lim_{t\uparrow 1}
\mu(x,t)dx =d\mu_1(x)$.
As a consequence,

\begin{eqnarray*}
\inf\{\alpha^*(\mu,M)+\beta^*(\mu,M)\}
&=&
\inf\{ \Sa(\rho,m)\quad : \quad (\rho,m)\mbox{ satisfies (C') and }
\rho_t|_{ K^c}=0\quad\forall t\in [0,1]\}\\
&=&2\inf\{S_{\mu_0}(\nu) \quad :\quad \nu_1=\mu_1\}
+(\Sigma(\mu_0)-\Sigma(\mu_1)):= Z(\mu_0,\mu_1)\\
\end{eqnarray*}
where in the last line we have used Property \ref{nfS}
and 
Corollary \ref{pr}.a).

Observe that $\a,\beta$ are convex functions with
values in $]-\infty,\infty]$. Moreover, there is at least one point 
$(F,B)\in E$, namely $F=-1,\quad\Phi=0$
for which $\a$ is continuous for the uniform topology on $E$ and
$\beta$ finite (this is the reason 
why we need to work on a compact set $K$ instead
of $\R$). Thus, following \cite{brenier},
by the Fenchel-Rockafellar duality theorem
(see th\'eor\`eme 1.11 in \cite{brezis}), we have

\begin{eqnarray*}
\inf\{\alpha^*(\mu,M)+\beta^*(\mu,M),\quad (\mu,M)\in E'\}
&=&\sup\{-\alpha(-
F,-\Phi)-\beta(F,\Phi) :
(F,\Phi)\in E\}
\end{eqnarray*}
and the infimum is achieved. 
More precisely,
$$Z(\mu_0,\mu_1)=\sup
\{ \int_Q\partial_t\phi_t(x) \rho^*_t(x) dxdt+
\int_Q \partial_x\phi_t(x) m^*_t(x)dx dt-{2\over 3\pi}
\int_Q \left(\partial_t\phi_t(x) +{1\over 4}\partial_x\phi_t(x)^2
\right)^{3\over 2}dxdt
\}$$
where the supremum 
is taken over  ${\phi\in\Ca^{1,1}_b(K\ts[0,1])}$ such that
$\partial_t\phi+{1\over 4}(\partial_x\phi)^2\ge 0.$
As a consequence, there exists a sequence of functions $\phi^\e$
in $\Ca_b^{1,1}(K\ts [0,1])$ such that 
$$\partial_t\phi^\e+{1\over 4}(\partial_x\phi^\e)^2\ge 0$$
and 
$$\int u^*_t(x)^2 d\mu^*_t(x)dt +{\pi^2\over 3}
\int \rho^*_t(x)^3 dxdt
\le \int \left( \partial_t\phi^\e_t(x)+u^*_t(x) 
\partial_x\phi^\e\right)d\mu^*_t(x)dt-{2\over 3\pi}
\int \left(\partial_t\phi_t^\e(x) +{1\over
4}\partial_x\phi_t^\e(x)^2
\right)^{3\over 2} +\e^2$$
for all $\e>0$, which implies, if we set
$$\pi^2(\rho^\e)^2=\partial_t\phi^e+4^{-1}(\partial_x\phi)^2,$$

\begin{eqnarray}
 \int (u^*_t(x)-\partial_x{\phi^\e_t(x)\over 2}
)^2 d\mu^*_t(x)dt
&\le& \pi^2
\int \rho^\e_t(x)^2 \rho^*_t(x)dxdt-{2\pi^2 \over 3}\int 
\rho_t^\e(x)^3 dx dt -{\pi^2\over 3}\int \rho^*_t(x)^3 dxdt
+\e^2\nonumber\\
&=&-{\pi^2\over 3}\int ( \rho^*_t(x)-\rho^\e_t(x))^2
(
2\rho^\e_t(x)+\rho^*_t(x))dxdt+\e^2\label{ty}\\
\nonumber
\end{eqnarray}
which completes the proof of the Property.

\hfill\xx

\section{Applications to matrix integrals}
In physics, several matrix integrals 
have been of interests in the 80's and 90's
for their applications to quantum fields
theory as well as string theory.
We refer here to the works
of M. Mehta, A. Matytsin, 
A. Migdal, V. Kazakov, P. Zinn Justin and B. Eynard for instance.
Among these integrals,
are often considered the following :

\begin{itemize}

\item The random Ising model on random
graphs described by the Gibbs measure

$$\mu^N_{Ising}(dA,dB)=
{1\over Z^N_{Ising}}
e^{N\tr(AB)-N\tr(P_1(A))-N\tr(P_2(B))}dA dB$$
with $Z^N_{Ising}$ the partition function
$$Z^N_{Ising}=\int e^{N\tr(AB)-N\tr(P_1(A))-N\tr(P_2(B))}dA dB$$
and two polynomial functions $P_1,P_2$.
The limiting free energy
for this model was calculated by M. Mehta \cite{mehta}
in the case $P_1(x)=P_2(x)= x^2+gx^4$
and integration holds over $\Ha_N$.
However, the limiting spectral measures
of $A$ and $B$ under $\mu^N_{Ising}$ were not considered in that paper.
A discussion about this problem
can be found
in P. Zinn Justin \cite{zinnjustin1}.

\item One can also define
the Potts model 
on random
graphs described by the Gibbs measure

$$\mu^N_{Potts}(dA_1,...,dA_q)=
{1\over Z^N_{Potts}}
\prod_{i=2}^q e^{N\tr(A_1A_i)-N\tr(P_i(A_i))}dA_i e^{ -N\tr(P_1(A_1))}
dA_1.$$
The limiting spectral measures 
of $(A_1,\cdots,A_q)$ are discussed 
in \cite{zinnjustin1} when $P_i=gx^3-x^2$ (!).

\item As a straightforward generalization, one can consider
matrices coupled by a chain following S. Chadha, G. Mahoux 
and M. Mehta \cite{mehta2}
given by
$$\mu^N_{chain}(dA_1,...,dA_q)=
{1\over Z^N_{chain}}
\prod_{i=2}^q e^{N\tr(A_{i-1}A_i)-N\tr(P_i(A_i))}dA_i e^{ -N\tr(P_1(A_1))}
dA_1.$$
$q$ can eventually go to infinity
as in \cite{MZ}.

\item Finally, we can mention the so-called 
induced QCD studied in \cite{matytsin}.
It is described, if $\L=[-q,q]^D\subset\Z^D$, by

$$\mu^N_{QCD}(dA_i,i\in\L)=
{1\over Z^N_{QCD}} \prod_{i\in\L}   \int 
e^{N\sum_{j=1}^{2D} \tr(U_j A_{i+e_j  }U_j^* A_i)} 
\prod_{j=1}^{2D}  dm^\beta_N(U_j )
\prod_{i\in\L} e^{ -N\tr(P(A_i))}dA_i $$
where $(e_j)_{1\le j\le 2D}$ is a basis of $\Z^D$. 
The description of the limit behaviour of the 
 spectral measures of $A_1,\cdots,A_q$ 
is given in \cite{matytsin} in the case $q=\infty$. We impose 
periodic boundary conditions 
at the boundary of the lattice points $\L$.

\end{itemize} 
In this section,
we shall study the asymptotic
behaviour of the free energy of these models 
as well as describe
the limit behaviour of 
the spectral measures of the matrices
under the corresponding Gibbs measures.

The theorem states as follows 

\begin{theo}\label{lim}
Assume that $P_i(x)\ge c_ix^4+d_i$ with $c_i>0$ and some finite
constants $d_i$. Hereafter, $\beta=1$ (resp. $\beta=2$) 
when $dA$ denotes 
the Lebesgue measure on $\Sa_N$ (resp. $\Ha_N$). Then,

\begin{eqnarray}
F_{Ising}&=& \lim_{N\ra\infty}
{1\over N^2}\log
Z^N_{Ising}\nonumber\\
&=&-\inf\{ \mu(P)+\nu(Q) -I^{(\beta)}(\mu,\nu)
-{\beta\over 2}
\Sigma(\mu)-{\beta\over 2}
\Sigma(\nu)\} -2\inf_{\nu\in\Pa(\R)} I_\beta
(\nu) \label{fIsing}\\
&&\nonumber
\\
F_{Potts}&=& \lim_{N\ra\infty}
{1\over N^2}\log
Z^N_{Potts}\nonumber\\
&=&-\inf\{ \sum_{i=1}^q \mu_i(P_i) -\sum_{i=2}^qI^{(\beta)}(\mu_1,\mu_i)
-{\beta\over 2}\sum_{i=1}^q 
\Sigma(\mu_i)\} -q\inf_{\nu\in\Pa(\R)} I_\beta
(\nu)\label{fPotts}\\
&&\nonumber
\\
F_{chain}&=& \lim_{N\ra\infty}
{1\over N^2}\log
Z^N_{chain}\nonumber\\
&=&-\inf\{ \sum_{i=1}^q \mu_i(P_i) -\sum_{i=2}^qI^{(\beta)}(\mu_{i-1},\mu_i)
-{\beta\over 2}\sum_{i=1}^q 
\Sigma(\mu_i)\} -q\inf_{\nu\in\Pa(\R)} I_\beta
(\nu) \label{fchain}\\
&&\nonumber
\\
F_{QCD}&=& \lim_{N\ra\infty}
{1\over N^2}\log
Z^N_{QCD}\nonumber\\
&=&-\inf\{ \sum_{i\in\L } \mu_i(P) -
\sum_{i\in\L}
\sum_{j=1}^{2D}
I^{(\beta)}(\mu_{i+e_j },\mu_i)
-{\beta\over 2}\sum_{i\in\L }
\Sigma(\mu_i)\} -2D|\L|\inf_{\nu\in\Pa(\R)} I_\beta
(\nu)\label{fqcd}\\
\nonumber
\end{eqnarray}
\end{theo}
\begin{rem} The above theorem actually extends 
to  polynomial functions going to infinity
like $x^2$. However, the case of quadratic 
polynomials is trivial since it boils down to
the Gaussian case and therefore the next interesting 
case is quartic
polynomial as above. Moreover, Theorem \ref{limitpoints} fails in the
case 
where $P,Q$ go to
infinity only like $x^2$.
However, all our proofs would extends easily for 
functions $P_i's$ such that $P_i(x)\ge a |x|^{2+\e}+b$
with some $a>0$ and $\e>0$.
\end{rem}

Theorem \ref{lim} will
be proved in the next section, but merely 
boils down to a Laplace's (or saddle point) method.

We shall then study 
the variational problems for the above energies.
We prove
the following for the Ising model.

\begin{theo}\label{limitpoints}\
Assume $P_1(x)\ge ax^4+b,P_2(x)\ge ax^4+b $ for some positive constant
$a$. Then
\begin{itemize}
\item The infimum in $F_{Ising}$
is achieved at a unique couple $(\mu_A,\mu_B)$ 
of probability measures.  

\item $(\mu_A,\mu_B)$ are compactly supported measures
with
finite entropy $\Sigma$.

\item Let $(\rho^{A\ra B}, u^{A\ra B})$ 
be the minimizer
of $S_{\mu_A}$ on $\{\nu_1=\mu_B\}$ 
as described in Theorem \ref{variational}. 
Then, $(\mu_A,\mu_B,\rho^{A\ra B},m^{A\ra B}=\rho^{A\ra B}u^{A\ra B})$
is the unique minimizer of the strictly convex  energy
\begin{eqnarray*}
\La(\mu,nu,\rho^*,m^*)&:=& \mu(P_1-{1\over 2} x^2)+
\nu(P_2-{1\over 2} x^2)
-{\beta\over 4}(
\Sigma(\mu)+
\Sigma(\nu))\nonumber\\&&
+{\beta\over 4}\left( \int_0^1\int {(m^*_t(x))^2\over \rho^*_t(x)}
 dxdt
+{\pi^2\over 3}\int_0^1\int \rho^*_t(x)^3 dx dt \right)\\
\end{eqnarray*}
Thus, we find that
$(\mu_A,\mu_B,\rho^{A\ra B},m^{A\ra B})$
are characterized by
the property  that for any $(\mu,\nu,\rho^*,m^*)\in\{\La<\infty\}$,

\begin{eqnarray}
&&\int (P_1-{1\over 2} x^2) d(\mu-\mu_A)+\int (P_2-{1\over 2} x^2)
 d(\nu-\mu_B)\nonumber\\
&&-{\beta\over 2} \int\int\log|x-y|d\mu_A(y)
(d\mu-d\mu_A)(x)
-{\beta\over 2} \int\int\log|x-y|d\mu_B(y)
(d\mu-d\mu_B)(x)\nonumber\\
&&+{\beta\over 4}
\int [2u^{A\ra B}(m^*-m^{A\ra B})-
(u^{A\ra B})^2 (\rho^*-\rho^{A\ra B})+
\pi^2 (\rho^{A\ra B})^2(\rho^*-\rho^{A\ra B})]dxdt\ge 0
\nonumber\\
\nonumber
\end{eqnarray}

\item  $(\rho^{A\ra B},m^{A\ra B})$
satisfies
 the Euler equation 
for isentropic flow with pressure $p(\rho)=-{\pi^2\over 3}\rho^3$
in the strong sense in the interior of  $\Omega=\{(x,t)\in\R\ts [0,1]: 
\rho^{A\ra B}_t(x)\neq 0\}$ and satisfy 
the conclusions of  Property \ref{variational}.

\item Moreover, 
$${1\over 2}\int\int {h(x)-h(y)\over x-y}d\mu_A(x)d\mu_A(y)
= \int (P_1'(x)-x -u_0^{A\ra B}(x))h(x)d\mu_A(x)$$
and
$${1\over 2}\int\int {h(x)-h(y)\over x-y}d\mu_B(x)d\mu_B(y)
= \int (P_2'(x)-x -u_1^{A\ra B}(x))h(x)d\mu_B(x)$$
for all $h\in\Ca^1_b(\R)$. 
\end{itemize}
\end{theo}
\begin{rem} The last point becomes
$$H\mu_A(x)=P_1'(x)-x -u_0^{A\ra B}(x),\quad\mu_A\mbox{ a.s.},
H\mu_B(x)=P_2'(x)-x -u_1^{A\ra B}(x),\quad\mu_B\mbox{ a.s.}$$
as soon as $H\mu_A$ and $H\mu_B$
are in $L^1(\mu_A)$ and  $L^1(\mu_B)$ respectively.
\end{rem}

For the other models, we unfortunately loose obvious
convexity, and therefore 
uniqueness of the minimizers in general.
We can still prove the following

\begin{theo}\label{limitpoints2}

\begin{itemize}

\item For any given $\mu_1$, there exists at most
one minimizer  $(\mu_2,\cdots,\mu_q)$ in $F_{Potts}$
but uniqueness of $\mu_1$ is  unclear in general, except in the case $q=3$
The critical points in $F_{Potts}$
are compactly supported, with finite entropy $\Sigma$.

Let $(\mu_1,\cdots,\mu_q)$ be a critical point
and for $i\in\{2,\cdots,q\}$, denote 
 $(\rho^i,u^i)$ the unique minimizer
described in Theorem \ref{variational}
with $\mu_0^i(dx)=\mu_1(dx)$
and $\mu_1^i(dx)=\mu_i(dx)$. Then 
$$P_1'(x)=qx +{\beta\over 2}\sum_{i=2}^q u^i_0 (x)
-{\beta\over 2}(q-3)H\mu_1(x)$$
in the sense of distributions on $\mbox{supp}(\mu_1)$
and
$$
P_i'(x)=x -{\beta\over 2}u^i_1(x)-{\beta\over 2}H\mu_i(x)
 ,\quad 2\le i\le q$$
in the sense of distributions on $\mbox{supp}(\mu_i)$.

\item There exists at most one minimizer in $F_{Chain}$,.
The minimizer $(\mu_1,\cdots,\mu_q)$ is compactly supported with
finite entropy $\Sigma$. 
The critical points    $(\mu_1,\cdots,\mu_q)$  in $F_{chain}$
are such that for $i\in\{2,\cdots,q\}$,W
It is such that  if we denote 
 $(\rho^i,u^i)$ the minimizer described in 
Theorem \ref{variational} 
with $\mu_0^i(dx)=\mu_{i-1}(dx)$
and $ \mu_1^i(dx) =\mu_i(dx)$,  we have 

$$P_1'(x)=x+{\beta\over 2}u^2_0-{\beta\over 2}H\mu_1(x)
\mbox{ and }
P_i'(x)=2x -{\beta\over 2}(u^i_1-u^{i+1}_0) , \quad 2\le i\le q.$$
in the sense of distribution in $\mbox{supp}(\mu_1)$ and $\mbox{supp}(\mu_i)$
respectively. 

\item Again, uniqueness
of the critical points in $F_{QCD}$
is unclear in general, except in the case $D=1$ where uniqueness holds.
In this case, the minimizer $\mu_i$ is symmetric, yielding
$\mu_i=\mu$ for all $i\in\L$ and 
 the unique path $(\rho^{},u^{})$ described
in  Theorem \ref{variational}  with  boundary data $(\mu,\mu)$,
satisfies $u^*_0(x)=-u^*_1(x)$ and 

$$P'(x)-2x
-\beta u_0^*(x)=0\quad \mu\mbox{a.s} $$

\end{itemize}
\end{theo}

\subsection{Proof of Theorem \ref{lim}}
The proof of Theorem \ref{lim}
follows a standard Laplace's method.
We shall only detail it in the Ising model
case, the generalization
to the other models being straightforward.

Let
$P,Q$ be two polynomial functions and define,
for $N\in\N$,
$\Lambda_N(P,Q)\in\R\cup\{+\infty\}$
by
$$\Lambda_N^\beta(P,Q)=\int \exp\{- N\tr(P(A))-N\tr(Q(B))
+N\tr(AB)\} dA dB$$
where the integration holds over orthogonal (resp. Hermitian)
 matrices if $\beta=1$
(resp. $\beta=2$).

We claim that

\begin{lem}
Assume that there exists $a,c\in\R^{+*}$, 
and $b,d\in\R$ such that
$$P(x)\ge ax^4+b\mbox{ and } Q(x)\ge cx^4+d,\qquad\mbox{ for all
}x\in\R.$$
Then,
we have 

$$\lim_{N\ra\infty}{1\over N^2}\log \Lambda_N(P,Q)
= \sup_{\mu,\nu\in \Pa(\R)}\{
-\mu(P)-\nu(Q) +I^{(\beta)}(\mu,\nu)
+{\beta\over
2}(\Sigma(\mu)+\Sigma(\nu) )\}-2\inf_{\nu\in \Pa(\R)} I_\beta(\nu)$$
\end{lem}
Remark here that the result could be extend to $
P(x)\ge ax^2+b\mbox{ and } Q(x)\ge cx^2+d$ with $ac>1$ but that
the Gaussian case being uninteresting, we shall use
the above and simpler hypothesis.

\nn
\prf

Observe that for any $\e>0$,
\begin{eqnarray*}
|\trn(AB)-\trn\left(
{A\over 1+\e A^2}{B\over 1+\e B^2}\right)|
&\le&\e \left|\trn\left( {A^3\over 1+\e A^2}
B\right)\right|+\e \left|\trn\left( {A\over 1+\e A^2} {B^3\over 1+\e B^2}
\right)\right|\\
&\le& \e\bigg(\left(\trn\left( {A^6\over (1+\e A^2)^2}
\right)\right)^{1\over 2} \left(\trn B^2\right)^{1\over 2}\\
&&
+ \left(\trn\left( {B^6\over (1+\e B^2)^2}
\right)\right)^{1\over 2} \left(\trn A^2\right)^{1\over 2}
\bigg)\\
&\le& \sqrt{\e}\left(\left(\trn(A^4)
\right)^{1\over 2} \left(\trn (B^2)\right)^{1\over 2}
+ \left(\trn( B^4)\right)^{1\over 2} \left(\trn A^2\right)^{1\over 2}
\right)\\
&\le& \sqrt{\e}\left(\trn(A^4)+\trn(B^4)+\trn(A^2) +\trn(B^2)
\right)\\
\end{eqnarray*}
Therefore,
if we set
$$\mu^N_{Ising}(dA,dB)={1\over\Lambda_N^\beta(P,Q)}
 \exp\{- N\tr(P(A))-N\tr(Q(B))+N\tr(AB)\}
dAdB$$
and
\begin{eqnarray*}
\D_N(\e)&:=&\left|{1\over N^2}\log
{
\int \exp\{- N\tr(P(A))-N\tr(Q(B))
+N\tr({A\over 1+\e A^2} {B\over 1+\e B^2})\} dA dB
\over \Lambda_N^\beta(P,Q)}\right|,\\
&=& \left|{1\over N^2}\log\mu^N_{Ising}
\left(  \exp\{N\tr({A\over 1+\e A^2} {B\over 1+\e
B^2})-N\tr(AB)\}\right)\right|,\\
\end{eqnarray*}
we get
\begin{eqnarray*}
\D_N(\e)&\le& {1\over N^2}\log \mu^N_{Ising} 
\left( \exp\{\sqrt{\e}N\tr( A^4+A^2)
+\sqrt{\e}N\tr( B^4+B^2)\}\right)\\
&\le & {1\over qN^2}\log \mu^N_{Ising} 
\left( \exp\{q\sqrt{\e}N\tr( A^4+A^2)
+q\sqrt{\e}N\tr( B^4+B^2)\}\right)\\
\end{eqnarray*}
where we  used Jensen's 
 inequality with $q>1$.
Now, under our hypothesis, 
and since  $2|AB|\le A^2+B^2$, it is clear that  if $q\sqrt{\e}$
is chosen small enough (e.g smaller than $a\wedge c$),
the above right hand side is bounded uniformly.
Hence, we take $q={1\over 2 a\wedge c\sqrt{\e}}$
and obtain

\begin{eqnarray}
\limsup_{N\ra\infty}
\D_N(\e)&\le& C\sqrt{\e}\label{delta}\\
\nonumber
\end{eqnarray}
with a finite constant $C$. 
Moreover, for any $\e>0$, we can use saddle point method
(see \cite{BAG} for a full rigorous derivation)  and 
Theorem 1.1 of \cite{GZ2} to obtain

$$\lim_{N\ra\infty}{1\over N^2}\log
\int \exp\{- N\tr(P(A))-N\tr(Q(B))
-N\tr({A\over 1+\e A^2} {B\over 1+\e B^2})\} dA dB$$
$$
= \sup_{\mu,\nu\in \Pa(\R)}\{
-\mu(P)-\nu(Q) +I^{(\beta)}(\mu\circ \phi_\e^{-1},\nu\circ \phi_\e^{-1})
+{\beta\over
2}(\Sigma(\mu)+\Sigma(\nu))\}-2\inf I_\beta$$
with $\phi_\e(x)=(1+\e x^2)^{-1} x$ and 
$\mu\circ\phi_\e^{-1} (f)=\mu(f\circ\phi_\e).$
Thus, (\ref{delta}) results with 

$$\lim_{N\ra\infty}{1\over N^2}\log \Lambda_N(P,Q)
=\lim_{\e\ra 0}
 \sup_{\mu,\nu\in \Pa(\R)}\{
-\mu(P)-\nu(Q) +I^{(\beta)}(\mu\circ \phi_\e^{-1},
\nu\circ \phi_\e^{-1})
+{\beta\over
2}(\Sigma(\mu)+\Sigma(\nu))\}-2\inf I_\beta.$$
Moreover, we can prove as for  (\ref{delta})
that
for any $\mu,\nu$ such that $\mu(x^4)\le M$
and $\nu(x^4)\le M$,
$$|I^{(\beta)}(\mu\circ \phi_\e^{-1},
\nu\circ \phi_\e^{-1})-I^{(\beta)}(\mu,
\nu)|\le C(M) \sqrt{\e}.$$
Using the fact that
$$|I^{(\beta)}(\mu,
\nu)|\le {1\over 2}(\mu(x^2)+\nu(x^2)),$$
as well as
$$\Sigma(\mu)+\Sigma(\nu)\le C(\mu(x^2)+\nu(x^2)+1)$$
for some finite constant $C$, 
we see that the supremum above is taken 
at $\mu,\nu$ such that $\mu(x^4)$ and $\nu(x^4)$
are bounded by some finite constant
depending only on $P,Q$. Hence, we can take the limit 
 $\e$ going to zero above and conclude.
\hfill\xx

\subsection{Proof of Theorem \ref{limitpoints} and \ref{limitpoints2}}
\subsubsection{The Ising model}
Let us recall that
\begin{eqnarray}
F_{Ising}+2\inf_{\nu\in\Pa(\R)} I_\beta 
(\nu)
&=&-\inf\{ \mu(P_1)+\nu(P_2) -I^{(\beta)}(\mu,\nu)
-{\beta\over 2}
\Sigma(\mu)-{\beta\over 2}
\Sigma(\nu)\}  \nonumber\\
\nonumber
\end{eqnarray}
Observe that since $I^{(\beta)}(\mu,\nu)\le
2^{-1}\mu(x^2)+2^{-1}\nu(x^2)$, 
the minimizer $(\mu_A,\mu_B)$
in the above right hand side is such that
$$\mu_A(P_1-{1\over 2} x^2)+\mu_B(P_2-{1\over 2} x^2)
-{\beta\over 2}
\Sigma(\mu_A)-{\beta\over 2}
\Sigma(\mu_B) \le -F_{Ising}-2\inf_{\nu\in\Pa(\R)} I_\beta 
(\nu)<\infty.$$
Hence, since $P_1-2^{-1} x^2$ and $P_2-2^{-1} x^2$
are bounded below under our hypotheses (for
well chosen $a$), we conclude 
that
$\Sigma(\mu_A)$ and $\Sigma(\mu_B)$ are bounded below
and hence finite. Further, if $2n_1$ (resp. $2n_2$)
is the degree  of $P_1$ (resp. $P_2$) for $n_1,n_2\ge 2$,
 we also see that
\begin{equation}\label{co1}
\mu_A(x^{2n_1})<\infty,\quad \mu_B(x^{2n_1})<\infty.
\end{equation}
Thus, we can use Property \ref{nfS}  to
get

\begin{eqnarray}
F_{Ising}+2\inf_{\nu\in\Pa(\R)} I_\beta 
(\nu)
&=& -\inf\bigg\lbc \mu(P_1-{1\over 2} x^2)+\nu(P_2-{1\over 2} x^2) 
-{\beta\over 4}(
\Sigma(\mu)+
\Sigma(\nu)) \nonumber\\
&&+
{\beta\over 4}\inf_{(u^*,\mu^*)\in (C)_{\mu,\nu}}
\{ \int_0^1\int u^*_t(x)^2 d\mu^*_t (x)dt
+\int_0^1\int H\mu^*_t(x)^2 d\mu^*_t(x) dt \} \bigg\rbc
 \nonumber\\
&=& -\inf_{\mu,\nu\in\Pa(\R)\atop (u^*,\mu^*)\in (C)_{\mu,\nu}}
\bigg\lbc \mu(P_1-{1\over 2} x^2)+\nu(P_2-{1\over 2} x^2)
-{\beta\over 4}(
\Sigma(\mu)+
\Sigma(\nu))\nonumber\\&&
+{\beta\over 4}\left( \int_0^1\int u^*_t(x)^2 d\mu^*_t (x)dt
+{\pi^2\over 3}\int_0^1\int \rho^*_t(x)^3 dx dt \right)\bigg)\nonumber\\
&:=& -\inf_{\mu,\nu\in\Pa(\R)\atop (u^*,\mu^*)\in (C)_{\mu,\nu}}
L(\mu,\nu,\mu^*,u^*)
\nonumber
\end{eqnarray}
where $(u^*,\mu^*)\in (C)_{\mu,\nu}$ means that
in the sense of distributions 
$$\partial_t\rho^*_t+\partial_x(\rho^*_t u^*_t)=0,\quad
\lim_{t\downarrow 0} \mu^*_t(dx) =\mu,\quad \lim_{t\uparrow 1}
 \mu^*_t(dx)=\nu$$
and we have used in the last line 
that when the above infimum is finite, $\mu^*_t$
is absolutely continuous with respect to
Lebesgue measure for almost all $t\in [0,1]$ and
with density $\rho^*\in L_3(dxdt)$ (see Lemma \ref{l3}).

\bigskip

Observe that if $L(\mu,\nu,\mu^*,u^*)=\La(\mu,\nu,\rho^*,m^*)$
with $m^*=\rho^*u^*$, $\La$ is a strictly convexe 
function of $(\mu,\nu, \rho^*,m^*)$ (recall that
$-\Sigma$ is convex, see \cite{BAG} for instance)
and that the constraint 
$(C)_{\mu,\nu}$ is linear in the variables $
(\mu,\nu,\rho^*,m^*)$.
Therefore, the above minimum
is achieved at a unique point 
$(\mu_A,\mu_B,\mu^{A\ra B}_., m^{A\ra B}_.)$.

\bigskip

We now perform a measure type perturbation to characterize the
infimum. Take
$ (\mu,\nu,\rho^*,m^*)\in \{\La<\infty\}$
and set, for $\a\in [0,1]$,
$$(\mu^\a,\nu^\a,\rho^\a,m^\a)=
\a (\mu,\nu,\rho^*,m^*)+(1-\a)(\mu_A,
\mu_B,\rho^{A\ra B}_., u^{A\ra B}_.).$$
Then, we find that we must have

\begin{eqnarray}
&&\int (P_1-{1\over 2} x^2) d(\mu-\mu_A)+\int (P_2-{1\over 2} x^2)
 d(\nu-\mu_B)\nonumber\\
&&-{\beta\over 2} \int\int\log|x-y|d\mu_A(y)
(d\mu-d\mu_A)(x)
-{\beta\over 2} \int\int\log|x-y|d\mu_B(y)
(d\mu-d\mu_B)(x)\nonumber\\
&&+{\beta\over 4}
\int [2u^{A\ra B}(m^*-m^{A\ra B})-
(u^{A\ra B})^2 (\rho^*-\rho^{A\ra B})+
\pi^2 (\rho^{A\ra B})^2(\rho^*-\rho^{A\ra B})]dxdt\ge 0
\label{eqising}\\
\nonumber
\end{eqnarray}
Taking $\mu=\mu_A$ and $\nu=\mu_B$, we see that
$(\rho^{A\ra B},u^{A\ra B})$
must satisfy Property \ref{variational}.
Now, if $\mu(dx)=\mu_A(dx)+\partial_x \phi_0(x) dx$,
 $\nu(dx)=\mu_B(dx)+\partial_x \phi_1(x) dx$
and $m^*=m^{A\ra B} -\partial_t\phi$, $\rho^*_t=
\rho^{A\ra B} +\partial_x\phi$ with $\phi\in\Ca^{\infty,\infty}_b(
\R\ts [0,1])$ such that
\begin{equation}\label{c12}
\int_{\rho^{A\ra B}\neq 0}{(m^{A\ra B}_t+\e\partial_t\phi)^2\over 
\rho^{A\ra B}+\e\partial_x\phi} dx dt<\infty,\quad
\int_{\rho^{A\ra B}= 0} {(\partial_t \phi)^2\over \partial_x\phi}
 dxdt<\infty,
\end{equation} 
we obtain 
by (\ref{eqising})
\begin{eqnarray}
&&\int (P_1-{1\over 2} x^2) \partial_x\phi_0(x)dx
+\int (P_2-{1\over 2} x^2)
  \partial_x\phi_1(x)dx\nonumber\\
&&-{\beta\over 2} \int\int\log|x-y|d\mu_A(y)
\partial_x\phi_0(x)dx-{\beta\over 2} \int\int\log|x-y|d\mu_B(y)
\partial_x\phi_1(x)dx\nonumber\\
&&+{\beta\over 4}
\int [2u^{A\ra B}\partial_t\phi -
(u^{A\ra B})^2 \partial_x\phi +
\pi^2 (\rho^{A\ra B})^2\partial_x\phi]dxdt\ge 0
\label{eqising2}\\
\nonumber
\end{eqnarray}
which becomes an equality  if $\phi$ is supported in $\Omega=\{(x,t)\in\R\ts [0,1] : \rho^{A\ra B}\neq 0\}$
by symmetry. If we assume that 
$u^{A\ra B}$ is sufficiently smooth, in particular 
 continuously differentiable with respect to
the time variable around $t=0$ and $t=1$, 
we can use integration by parts to see
that
$$
\int [2u^{A\ra B}\partial_t\phi -
(u^{A\ra B})^2 \partial_x\phi +
\pi^2 (\rho^{A\ra B})^2\partial_x\phi]dxdt\ge 2[\int \Pi^{A\ra B}_t \partial_x\phi_t dx]_0^1
$$
yielding that there exists two  constants $l_1,l_2$
 such that 
\begin{eqnarray}
P_1(x)-{1\over 2} x^2-
{\beta\over 2}\int\log|x-y|d\mu_A(y)-
2\Pi^{A\ra B}_0(x)&=&l_1\quad \mu_A\quad a.s
\label{f1}\\
P_2(x)-{1\over 2} x^2-{\beta\over 2}\int\log|x-y|d\mu_B(y)
-2\Pi^{A\ra B}_1(x)&=&l_2\quad \mu_B\quad a.s\label{f2}\\
P_1(x)-{1\over 2} x^2-{\beta\over 2}\int\log|x-y|d\mu_A(y)
-2\Pi^{A\ra B}_0(x)&\ge&l_1\quad \mbox{ if }x\in\mbox{supp}(\mu_A)^c
\nonumber\\
P_2(x)-{1\over 2} x^2-{\beta\over 2}
\int\log|x-y|d\mu_B(y)-2\Pi^{A\ra B}_1(x)&\ge&l_2\quad
 \mbox{ if }x\in\mbox{supp}(\mu_B)^c\nonumber\\
\nonumber
\end{eqnarray}
Such a result would generalize the usual equations obtained in 
the one matrix case. However, since we
could not prove  such a regularity 
property of $(\rho^{A\ra B},u^{A\ra B})$,
we shall now obtain a 
Schwinger-Dyson type 
formula following \cite{CDG2}, theorem 2.15 and proposition 2.17,
to obtain a weak form of (\ref{f1}),(\ref{f2}).
Let us briefly recall the ideas in the case $\beta=2$
(the case $\beta=1$ being similar),
which is based on an infinitesimal change of variables.

If, in $Z^N_{Ising}$,  we change $A\ra A+N^{-1} h(A,B)$
 with some 
smooth bounded  functions $h$ of two non-commutative variables
(take for instance $h$ belonging
to the set $\Ca\Ca_{st}(\C)$ of Stieljes functionals  defined in
\cite{CDG1,CDG2}(see also its definition in appendice \ref{app1}),
it turns out that, due to Kadison-Fuglede determinant formula (see 
\cite{CDG2}, the proof of theorem 2.15 and proposition 2.17)

$$Z^N_{Ising}=\int e^{ \tr( h(A,B)(-P_1'(A)+B))
+N^{-1}\tr\otimes \tr (D_A h(A,B))+O(1)-N\tr(P_1(A)+P_2(B)-AB)} dAdB
$$
with
$D_A$ the non commutative derivation 
with respect to $A$ given by 
$$D_A (hg)=D_A h\ts 1\otimes g+h\otimes 1\ts D_A g,\quad\forall
h,g\in \Ca\Ca_{st}(\C), \quad
D_A B=0, \quad D_AA=1\otimes 1.$$ 
Therefore,
we can find a finite constant $C(h)$ such that
for any $\e>0$

\begin{equation}\label{cv}
\mu^N_{Ising}\left( 
|\hat\mu^{(N)}\otimes \hat\mu^{(N)}(D_A h(A,B))
+\hat\mu^{(N)}((-P_1'(A)+B) h(A,B))|\le \e\right)
\ge 1-2e^{-\e N+C(h)}
\end{equation}
with $ \hat\mu^{(N)}$
the empirical distribution of $A,B$ defined by
$$\hat\mu^{(N)}(h)=\trn(h(A,B)),\quad \forall h\in \Ca\Ca_{st}(\C).$$
Of course, the same type of formula holds 
when $A$ is replaced by $B$.  It is not hard to
see that $\hat\mu^{(N)}$ is tight under
$\mu^N_{Ising}$ for the topology described in 
\cite{CDG2}, corresponding to the $ \Ca\Ca_{st}(\C)$-weak
 topology (see \cite{CDG2} for proof of
similar statements).  Let $\tau$ be a limit point.
Taking, for $\e>0$ and $\d>0$,
$h(A,B)= (1+\d A^2)^{-p} j(A) (1+\e B^2)^{-1}$
with $j(x)=\prod_{1\le i\le n}(z_i-x)^{-1}$ for some
$z_i\in\C\backslash\R$ and $n\in\N$,
and $p$ large enough ($p$ larger than half the degree of $P_1'$)
so that $D_A h(A,B)\in \Ca\Ca_{st}(\C)\otimes\Ca\Ca_{st}(\C)$ 
and $(1+\d A^2)^{-p}(P_1'(A)-B) (1+\e B^2)^{-1} j(A)
\in \Ca\Ca_{st}(\C)$, we deduce from (\ref{cv})
that 
$\tau$ must satisfies for any  $\e,\d>0$
and $p$ large enough,
\begin{equation}\label{cv2}
\tau\otimes\tau( D_A (1+\d A^2)^{-p} j(A)\ts 1\otimes (1+\e
B^2)^{-1})
= \tau( (P_1'(A)-B) (1+\d A^2)^{-p} j(A) (1+\e B^2)^{-1}).
\end{equation}
Similarly for any $\e,\d>0$, 
and $p$ large enough,
\begin{equation}\label{cv3}
\tau\otimes\tau( D_B (1+\d B^2)^{-p}j(B) \ts 1\otimes (1+\e
B^2)^{-1})
= \tau( (P_2'(B)-A) (1+\d A^2)^{-p} j(B) (1+\e A^2)^{-1}).
\end{equation}
Now, by (\ref{co1}), $P_1'(A)-B$ and $P_2'(B)-A$ belongs
to $L^1(\tau)$ so that we can let $\d,\e$ going
to zero to
conclude by dominated convergence theorem that

\begin{equation}\label{cv4}
\tau\otimes\tau( D_A  j(A))
= \tau( (P_1'(A)-B)  j(A)),\quad \tau\otimes\tau( D_B j(B))
= \tau( (P_2'(B)-A) j(B) ).
\end{equation}
We next show that (\ref{cv4}) implies that $\mu_A$ and $\mu_B$ 
are compactly supported when
$n_1\ge 2$ and $n_2\ge 2$,
 and first that all their moments are finite.
To this end, take $j(x)=\left( (1+\e x^2)^{-1} x\right)^n
=\e^{-n}\left(1+i\e^{-1}(x-i\e^{-1})^{-1}
\right)^n (x+i\e^{-1})^{-n}$ for $n\in\N$, yielding
\begin{equation}\label{eqa1}
\mu_A\left( P_1'(x) \left( (1+\e x^2)^{-1} x\right)^n
\right) =\tau \left (\tau(B|A) 
\left( (1+\e A^2)^{-1} A\right)^n\right) + \tau\otimes\tau( D_A 
j(A))
\end{equation}
with, since $Df$ can  be represented in the
tensor product space as
$Df(x,y)=(x-y)^{-1}(f(x)-f(y))$,
$$\tau\otimes\tau( D_A  j(A))
=\sum_{p=0}^{n-1}\mu_A \left( ((1+\e x^2)^{-1} x)^p\right)
\mu_A\left( ((1+\e x^2)^{-1} x)^{n-1-p}\right)\qquad\qquad$$
$$\qquad\qquad
-\e \sum_{p=0}^{n-1} \mu_A \left( ((1+\e x^2)^{-1} x)^{p+1}\right)
\mu_A\left( ((1+\e x^2)^{-1} x)^{n-p}\right).$$
When $n$ is odd, it is not hard to
see that
we can find $c>0,d_n\in\R$ such that
$P'(x) x^n\ge cx^{2n_1-1+n}-d_n$, so that
we deduce from (\ref{eqa1})
that
\begin{equation}\label{eqa2}
c\mu_A\left( |{x\over 1+\e x^2}|^{2n_1-1+n}
\right)\le d_n+2n
\sup_{p\le n} \mu_A\left( ({x\over 1+\e x^2})^p\right)^2
+\mu_A\left( | {x\over 1+\e x^2}
|^{nq}\right)^{1\over q}
\mu_B(|x|^p)^{1\over p}
\end{equation}
where we have used in the last line H\"older's inequality with
conjugate exponents $p,q$. We take
$q=n^{-1}(2n_1-1+n)$, $p=(2n_1-1)^{-1}(2n_1-1+n)$.
Similarly, we obtain for $\mu_B$, and $q=n^{-1}(2n_2-1+n)$, 
$p=(2n_2-1)^{-1}(2n_2-1+n)$,
\begin{equation}\label{eqa3}
c\mu_B\left( |{x\over 1+\e x^2}|^{2n_2-1+n}
\right)\le d_n+2n
\sup_{p\le n} \mu_B\left( ({x\over 1+\e x^2})^p\right)^2
+\mu_B\left( |{x\over 1+\e x^2}|^{nq }\right)^{1\over 
q}
\mu_A(|x|^{p})^{1\over p}.
\end{equation}
Now, we have seen that
$$\mu_A(x^{2n_1})<\infty,\quad \mu_B(x^{2n_2})<\infty$$
so that (\ref{eqa2}),(\ref{eqa3}) yields 
\begin{eqnarray*}
\mu_A(x^{2n_1-1+n})&=& \sup_{\e\ge 0} \mu_A\left( ((1+\e x^2)^{-1} x)^{2n_1-1+n}
\right)<\infty \mbox{ for } {2n_1-1+n}\le m_1^A:= 2n_2(2n_1-1)\\
\mu_B(x^{2n_2-1+k})&=& \sup_{\e\ge 0} \mu_B\left( ((1+\e x^2)^{-1}
 x)^{2n_2-1+k}
\right)<\infty \mbox{ for } {2n_2-1+k}\le m_1^B:=2n_2(2n_2-1)\\
\end{eqnarray*}
and then by induction for $2n_1-1+n\le m_p^A:=m_{p-1}^B(2n_1-1)$,
${2n_2-1+k}\le m_p^B:=m_{p-1}^A (2n_2-1)$ for
all $p\ge 2$.
Since $ 2n_1-1> 1$ and $2n_1-1>1$, $m_p^A$ and $m_p^B$
go
to infinity with $p$, which proves that
$\mu_A$ and $\mu_B$ have finite moments of all orders.

\bigskip

As a consequence, we can extend 
by dominated convergence theorem (\ref{eqa1})
to polynomial functions (i.e. take $\e=0$)
resulting with

\begin{equation}\label{eqa4}
\mu_A\left( P_1'(x)  x^n
\right) =\tau \left (\tau(B|A) 
 A^n\right) +
\sum_{p=0}^{n-1}\mu_A \left(  x^p\right)\mu_A
\left( x^{n-1-p}\right).
\end{equation}
and a similar equation for the moments of $\mu_B$.
Let us write $P_1'(x)=\a_1 x^{2n_1-1}
+\sum_{p=2}^{2n_1} \a_p x^{2n_1-p}$, $P_2'(x)
=\b_1 x^{2n_2-1}
+\sum_{p=2}^{2n_2} \b_p x^{2n_2-p}$ with $\a_1>0,\b_1>0$.
Setting $a_n=|\mu_A(x^n)|$ and $b_n=|\mu_B(x^n)|$,
we deduce that

\begin{eqnarray}
\a_1 a_{2n_1-1+n}
&\le &
\sum_{p=2}^{2n_1}|\a_p|a_{2n_1-p+n}
+\sum_{p=0}^{n-1} a_p a_{n-1-p}
+ a_{qn}^{1\over q} b_p^{1\over p}\label{eqaz1}\\
\b_1 b_{2n_1-1+n}
&\le &
\sum_{p=2}^{2n_2}|\b_p|b_{2n_1-p+n}
+\sum_{p=0}^{n-1} b_p b_{n-1-p}
+ b_{qn}^{1\over q} a_p^{1\over p}\label{eqaz2}\\
\nonumber
\end{eqnarray}
with conjuguate exponents $(p,q)$ to
be chosen later. 
 
Now, we make the induction hypothesis that
for some $R\in\R^+$, for some $m\in\N$,
$$a_p\le R^p C_p,\quad b_p\le R^p C_p,\quad\mbox{ for }p\le m$$
with $C_p$ the Catalan numbers given by
$$C_p=\sum_{n=0}^{p-1} C_n C_{p-1-n},\quad C_0=1.$$
Of course, up to take $R$ big enough, we can always assume
that
$m\ge 2n_1\vee n_2$. 
Now, plugging this hypothesis into
(\ref{eqaz1}),(\ref{eqaz2}) with $m+1=2n_1-1+n$
and $q=mn^{-1}$, we obtain
\begin{eqnarray*}
\a_1 a_{2n_1-1+n}
&\le &\sum_{p=2}^{2n_1}|\a_p | 
R^{2n_1-p+n}C_{2n_1-p+n}
+R^{n}C_n
+ R^{n+1} (C_{m})^{n\over m} (C_{[{m\over m-n}]+1})^{m-n\over m}\\
&\le& C_{m+1} R^{m+1} (\sum_{p=2}^{2n_1} |\a_p|R^{-p} 
+R^{n-2-m} +R^{n-m})\\
\end{eqnarray*}
where we have used that $C_m$ increases with $m$.
Thus, our induction hypothesis is verified
as soon as
\begin{eqnarray*}
\sum_{p=2}^{2n_1} |\a_p|R^{-p} 
+R^{-2n_1} +R^{2(1-n_1)}&\le& \a_1\\
\sum_{p=2}^{2n_2} |\b_p|R^{-p} 
+R^{-2n_1} +R^{2(1-n_2)}&\le& \b_1\\
\end{eqnarray*}
which is clearly the case for $R$ large enough
since we asumed $n_1\wedge n_2\ge 2$. 
Since $m^{-1}\log C_m$ goes to $4$ as $m$ goes to infinity,
we deduce that
$$\limsup_{m\ra\infty}{1\over 2m}\log \mu_A(x^{2m})\le
R+4,\quad \limsup_{m\ra\infty}{1\over 2m}\log \mu_B(x^{2m})\le
R+4,$$
implying that $\mu_A$ and $\mu_B$ 
are supported into $[-R-4,R+4]$
for $R$ finite satisfying the above induction hypothesis
(plus the condition imposed by the first $2n_1\vee n_2$
moments). 

\bigskip

Let us now go back to (\ref{cv4}) and notice that
since the Stieljes functions are dense
in $\Ca_c(\R)$ and $P'_1-\tau(A|B)$ belongs to $L^1(\tau)$,
it can be extended to $j\in \Ca_b^1(\R)$ ;
\begin{equation}\label{cv6}
\int \int {j(x)-j(y)\over x-y} d\mu_A(x) d\mu_A(y)
=\tau ((P'(x)-\tau(B|A))j(x))
\end{equation}

\bigskip

Since
$(\mu_A,\mu_B)$
are compactly supported, 
we can use the conclusions of 
 section 
\ref{bbp}. We see that $\mu^{A\ra B}_t$-almost surely,

$$u_t^{A\ra B}= \tau( B-A|X_t) +(1-2t) H\mu_t^{A\ra B}(x)$$
so that 
$$u_0^{A\ra B}=\tau(B|A) -x+H\mu_A$$
at least in the sense of distribution as in (\ref{cv6}).
Thus, by uniqueness of the solutions to the Euler
equation given
the initial and final data $(\mu_A,\mu_B)$ proved in Property \ref{unique},
 we conclude that
$$H\mu_A(x)=P_1'(x)-x -u_0^{A\ra B}(x)$$
in the sense of distribution
that 
$${1\over 2}\int\int {h(x)-h(y)\over x-y}d\mu_A(x)d\mu_A(y)
= \int (P_1'(x)-x -u_0^{A\ra B}(x))h(x)d\mu_A(x)$$
for all $h\in\Ca^1_b(\R)$. 
The second equation is derived similarly
and one finds that 
$$2H\mu_B(x)=P_2'(x)-\tau(A|B)(x)=
P_2'(x)-(x-u_1^{A\ra B}(x)-H\mu_B(x))$$
resulting with
$$H\mu_B(x)=P_2'(x)-x+u_1^{A\ra B}(x)$$
in the sense of distributions. 
Note also that by Property \ref{variational}, the fact that
$(\mu_A,\mu_B)$ are compactly supported
implies that $(\rho^{A\ra B},u^{A\ra B})$
satisfies the isentropic Euler equation in the strong sense
in $\Omega$.

\subsection{q-Potts model}
In this case,
we find that

\begin{eqnarray}
F_{Potts}
&=&-\inf\{ \sum_{i=1}^q \mu_i(P_i) -\sum_{i=2}^qI^{(\beta)}(\mu_1,\mu_i)
-{\beta\over 2}\sum_{i=1}^q 
\Sigma(\mu_i)\} -q\inf_{\nu\in\Pa(\R)} I_\beta
(\nu)\nonumber\\
&=&-\inf\{ \mu_1(P_1-{q\over 2}x^2)+\sum_{i=2}^q \mu_i(P_i-{x^2\over 2})\nonumber \\
&&+{\beta\over 4} \sum_{i=2}^q
\inf_{(u^*,\mu^*)\in (C)_{\mu_1,\mu_i}}
\left\lbc \int_0^1\int (u^*_t(x))^2 d\mu^*_t (x)dt
+\int_0^1\int (H\mu^*_t(x))^2 d\mu^*_t(x) dt \right\rbc\nonumber\\
&&-{\beta\over 4}\sum_{i=2}^q 
\Sigma(\mu_i)+{\beta\over 4}(q-3) \Sigma(\mu_1)\} -q\inf_{\nu\in\Pa(\R)} I_\beta
(\nu)\label{fPotts2}\\
\nonumber
\end{eqnarray}
When $q> 3$, the above functionnal is not
 anymore clearly convex in $\mu_1$ since $\Sigma(\mu_1)$ is concave.
Hence, the uniqueness of the minimizers is
now unclear. Note however that the Euler-Lagrange 
term may still contain sufficient  convexity 
in $\mu_1$ to insure uniqueness but simply that 
the above formula does not
show it. In the case $q=3$, the functional is still convex, and strictly 
convex in the arguments $(\mu^*,m^*)$. Therefore, uniqueness
of the minimizers still holds since if $(\mu,\nu,\mu^*,u^*)$ 
and $(\tilde\mu,\tilde\nu,\tilde\mu^*,\tilde u^*)$, we would still 
find that by convexity $\tilde\mu^*=\mu^*$ and therefore 
$\mu=\mu^*_0=\tilde\mu^*_0=\tilde\mu$,
$\nu=\mu^*_1=\tilde\mu^*_1=\tilde\nu$.
The above formula
already shows that the critical 
points satisfy $\mu_i(P_i)<\infty$ 
and have finite entropy $\Sigma$.  We can also obtain the
Schwinger-Dyson equations
in this case and deduce as for the Ising model 
that the critical points  are compactly supported and
satisfy 
the equations of Theorem \ref{limitpoints2}.

\subsection{Chain model }
In this case,
\begin{eqnarray}
F_{chain}
&=&-\inf\{ \sum_{i=1}^q \mu_i(P_i) -\sum_{i=2}^{q} I^{(\beta)}(\mu_{i-1},\mu_i)
-{\beta\over 2}\sum_{i=1}^q 
\Sigma(\mu_i)\} -q\inf_{\nu\in\Pa(\R)} I_\beta
(\nu) \label{fchain2}\\
&=&-\inf\{\mu_1(P_1-{x^2\over 2})+ \sum_{i=2 }^q \mu_i(P_i-x^2)\nonumber\\
&&+{\beta\over 4} \sum_{i=2}^q
\inf_{(u^*,\mu^*)\in (C)_{\mu_i,\mu_{i+1}}}
\left\lbc \int_0^1\int (u^*_t(x))^2 d\mu^*_t (x)dt
+\int_0^1\int (H\mu^*_t(x))^2 d\mu^*_t(x) dt \right\rbc\nonumber\\
&&-{\beta\over 4}\Sigma(\mu_1)\} -q\inf_{\nu\in\Pa(\R)} I_\beta
(\nu) \label{fchain3}\\
\nonumber
\end{eqnarray}
Here, we still have convexity and strict convexity on the term
coming from $I^{(\beta)}$. Hence, uniqueness of the minimizers
hold. 
Again, we can prove the
conclusions of Theorem \ref{limitpoints2}
as for the Ising model.

\subsection{Induced QCD model}
\begin{eqnarray}
F_{QCD}
&=&-\inf\{ \sum_{i=1 }^q \mu_i(P) -
\sum_{i\in\L}
\sum_{j=1}^{2D}
I^{(\beta)}(\mu_{i+e_j},\mu_i)
-{\beta\over 2}\sum_{i\in\L }
\Sigma(\mu_i)\} -2D|\L|\inf_{\nu\in\Pa(\R)} I_\beta
(\nu)\nonumber\\
&=& -\inf\{ \sum_{i=1 }^q \mu_i(P-D x^2)-{\beta\over 2}(1-D)\sum_{i\in\L }
\Sigma(\mu_i)\nonumber\\
&& +
\sum_{i\in\L}
\sum_{j=1}^{2D}
\inf_{(u^{i,\mu},\mu^{i,\mu})\in (C)_{\mu_i,\mu_{i+\mu}}}
\left\lbc \int_0^1\int (u^{i,\mu}_t(x))^2 d\mu^{i,\mu}_t (x)dt
+\int_0^1\int (H\mu^{i,\mu}_t(x))^2 d\mu^{i,\mu}_t(x) dt \right\rbc\}\nonumber\\
&& -2D|\L|\inf_{\nu\in\Pa(\R)} I_\beta
(\nu)\nonumber\\
\nonumber
\end{eqnarray}
Again, obvious convexity disappears 
and uniqueness of the minimizers
becomes unclear when $D>1$.
Uniqueness of the minimizers
still holds when $D=1$. 
Then, clearly $\mu_i=\mu$ for all $i\in\L$
and  $u_0^*=-u_1^*$ at the minimizing path
with $(\rho^*,u^*)$ the solution of the Euler 
equation with with boundary data $(\mu,\mu)$.
$\mu$ then satisfies 
\begin{eqnarray*}
P'(x)-2x-{\beta}u_0^*(x) &=&0\\
\end{eqnarray*}
in the sense of distributions in $\mbox{supp}(\mu)$, 
which corresponds to the result obtained
 by Matytsin [\cite{matytsin}, (4.3)]
when $\beta=2$. Actually, since we can prove as 
for the Ising model that $\mu$ is compactly supported,
it turns out that $P'(x)-2x-{\beta}u_0^*(x)$ is in every $L_p(d\mu)$
and therefore that $P'(x)-2x-{\beta}u_0^*(x)=0$ almost everywhere
in the support of $\mu$.

\section{Appendice }

\subsection{ Free Brownian motion description
of the minimizers}\label{app1}

Let us return to the probability aspect of the story.
In fact, by definition, if 
$$X^N_t= X^N_0+ H^N_t$$
with a Hermitian (if $\b=2$, otherwise
symmetric if $\b=1$)
 matrix $ X^N_0$ with spectral measure 
$\mun_0$ and a Hermitian (resp. symmetric) 
Brownian motion
$H^N$, if we denote $\mun_t$ the spectral measure 
of $X^N_t$, then, if $\mun_0$ converges towards
a compactly supported probability measure
$\mu_0$,
for any $\mu_1\in\Pa(\R)$,
$$\limsup_{\d\ra 0}\limsup_{N\ra\infty} {1\over N^2}\log\P(d(\mun_1,\mu_1)<\d)
=\liminf_{\d\ra 0}\liminf_{N\ra\infty} {1\over N^2}\log\P(d(\mun_1,\mu_1)<\d)
=-J_\beta(\mu_0,\mu_1).$$
Let us now reconsider the
above limit and show that the
limit must be taken at a free Brownian bridge.
More precisely, we shall see that,
if $\tau$ denotes the joint law of $(X_0,X_1)$ (the
precise sense of which being given below)
and $\mu^\tau$ the law of 
the free Brownian bridge (\ref{fb})
associated with $(X_0,X_1)$,

$$\limsup_{\d\ra 0}\limsup_{N\ra\infty} {1\over N^2}\log\P(d(\mun_1,\mu_1)<\d)
\le \sup_{ \tau\circ X_0^{-1} =\mu_0\atop
\tau\circ X_1^{-1} 
=\mu_1} \limsup_{\d\ra 0}\limsup_{N\ra\infty} {1\over
N^2}\log\P(\max_{1\le k\le n}d(\mun_{t_k},\mu_{t_k}^\tau
)\le \d)$$
for any family $\{t_1,\cdots,t_n\}$ of times
in $[0,1]$. Therefore, the large deviation estimate 
obtained in \cite{GZ2} yields

$$\limsup_{\d\ra 0}\limsup_{N\ra\infty} {1\over N^2}\log\P(d(\mun_1,\mu_1)<\d)
\le -{\beta\over
2} \inf\{ S(\mu^\tau),\tau\circ X_0^{-1} =\mu_0,\tau\circ X_1^{-1} 
=\mu_1\}.$$
The lower bound estimate
obtained in \cite{GZ2}  therefore guarantees
that
$$ \inf\{ S(\nu),\nu_0=\mu_0,\nu_1=\mu_1\}=
 \inf\{ S(\mu^\tau),\tau\circ X_0^{-1} =\mu_0,\tau\circ X_1^{-1} 
=\mu_1\}.$$
Such kind of result were already obtained in \cite{CDG2}
and \cite{BCG}.

Let us now be more precise. We  recall that
we can define the joint law 
of the two matrices $X^N_0,X^N_1$
by the family
$$\mun_{0,1}(F)=\trn(F(X^N_0,X^N_1))$$
when $F$ is taken into
a natural set $\Fa$  of test functions
of two non-commutatives variables
and $\trn(A)=N^{-1}\sum_{i=1}^N A_{ii}$. It
is common in free probability to consider polynomial test functions.
In \cite{CDG1}, bounded analytic
test functions were introduced for self-adjoint  non-commutatives
 variables. $\Fa=\Ca\Ca_{st}(\C)$  is there
the complex vector space generated by
$$
F(X_1,X_2)=\prod_{1\le i\le n}^{\ra}
 \frac{1}{z_i-
\alpha_i^1X_1-\a_i^2X_2}
$$
where $(z_i)_{1\le i\le n}$ belongs to $(\C\backslash\R)^n$,
 $(\alpha_i^k,1\le k\le 2)_{i=1}^n$ to $(\R^2)^{n}$, and $\displaystyle\prod^{\ra}$ is 
the non-commutative product.

We shall here use the very same set of functions
and recall then that the space
$$\Ma_{0,1}=\{\tau\in\Fa^* :\tau(I)=1,\tau(FF^*)\ge 0, 
\tau(FG)=\tau(GF)\}$$
is a compact metric space.
We denote by $D$ a metric on $\Ma_{0,1}$.
Let us recall \cite{CDG1} that if one considers 
the restriction 
$\mu_k=\tau\circ X_{k}^{-1}$ of $\tau$ to functions 
which only depends on one of the variables $X_k$, $k=1,2$, 
then $\mu_k$ is a probability measure 
on $\R$ (in fact the spectral measure of $X_k$)
and that the topology inherited by duality from $\Fa$ is the 
{\it vague} topology, i.e. the topology 
generated by continuous compactly supported
functions.

Since $\Ma_{0,1}$ is compact, for any
$\e>0$, we can find $M\in\N$,  $(\tau_k)_{1\le k\le
M}$ so that $\Ma_{0,1}\subset
\cup_{1\le k\le M}\{\tau:D(\tau,\tau_k)<\e\}$
and therefore
\begin{eqnarray*}
\limsup_{\d\ra 0}\limsup_{N\ra\infty} {1\over N^2}\log\P(d(\mun_1,\mu_1)<\d)
&\le& \max_{1\le k\le M}
\limsup_{N\ra\infty} {1\over N^2}\log\P(d(\mun_1,\mu_1)<\d;
D(\mun_{0,1},\tau_k)<\e)\\
\end{eqnarray*} 
Now, conditionnally to $X^N_1$,

$$dX^N_t= dH^N_t -{X^N_t-X^N_1\over 1-t} dt$$
or equivalently

$$X^N_t=tX^N_1+(1-t)X^N_0 +(1-t)\int_0^t(1-s)^{-1} dH^N_s.$$
Let us assume that $\mun_{0,1}$
converges towards $\tau\in\Ma_{0,1}$ when $N$ goes to infinity
and 
 that $X^N_1,X^N_0$
remains uniformly bounded for the operator norm. In particular,
$\mun_{tX^N_1+(1-t)X^N_0}$ converges for any $t\in [0,1]$
towards $\nu_t^\tau=\tau\circ(tX_1+(1-t)X_0)^{-1}$;
$$\nu_t^\tau(f)=\tau(f(  tX_1+(1-t)X_0))$$
for any test function $f$.
Therefore, Voiculescu's result implies that
$\mun_{X^N_t}$ converges towards the distribution $\mu^\tau_t$
of $
tX_1+(1-t)X_0+(1-t)\int_0^t(1-s)^{-1} dS_s$
with a free Brownian motion $S$, free with $tX_1+(1-t)X_0$.
We shall now extend this result in
our topology and also 
control the dependence of this convergence with
respect to the 
speed of convergence of the distribution of $(X_0^N,X_1^N)$ 
towards $\tau$.

We shall work below with given $(X^N_0,X^N_1)\in
\{d(\mun_0,\mu_0)<\d; d(\mun_1,\mu_1)<\d;
D(\mun_{0,1},\tau)<\e\}$.

Let, for $u\le t$, $X^{N,t}_u$
denote the process
$$X^{N,t}_u= tX^N_1+(1-t)X^N_0 +(1-t)\int_0^u(1-s)^{-1} dH^N_s.$$
Then, one deduces from Ito's calculus  that
for any test function
$f$ 
$$\mun_{X^{N,t}_u}(f)=\mun_{tX^N_1+(1-t)X^N_1}(f)
+{(1-t)^2\over 2}\int_0^u \mun_{X^{N,t}_s}\otimes
\mun_{X^{N,t}_s}({f'(x)-f'(y)\over x-y}){ ds\over
(1-s)^2} +M^N_f(u)$$
with a martingale $M^N_f(u)$
such that 
$$\E[\sup_{u\in [0,t]} (M^N_f(u))^2]\le {||f'||_\infty^2\over N^2}.$$
Moreover, it is not hard to check that
$(\mun_{X^{N,t}_u},u\le t)$ is tight
in $\Ca([0,1],\Pa(\R))$ (see the proof of exponential tightness
of the spectral process of $X^N_0+H^N_t$
given in \cite{GZ2}). The limit points
$(\mu_{X^t_u}, u\le t)$ (when $D(\mun_{0,1},\tau)$ goes to zero)
 satisfy
the equation
$$\mu_{X^t_u}(f)=\nu_{t}^\tau(f)
+{(1-t)^2\over 2}\int_0^u \mu_{X^{t}_s}\otimes
\mu_{X^{t}_s}({f'(x)-f'(y)\over x-y}) {ds\over
(1-s)^2}.$$
This equation admits a unique solution,
as can be proved following
the arguments of \cite{CDG2} or \cite{GZ2}, p. 494.
Taking $f(x)=e^{i\xi x}$,
and substracting both equations,
we find, with
$$\D_u^N(R)=\sup_{|\xi|\le R}\E[|\mun_{X^{N,t}_u}(e^{i\xi x} )-
\mu_{X^t_u}(e^{i\xi x})|],$$
that for $u\le t$
$$\D_u^N(R)\le \D_0^N(R)+4R^2\int_0^u \D_s^N(R)ds
+{R\over N}$$
which yields thanks to Gronwall lemma and taken at $u=t$,
since  $\mu^\tau_t=\mu_{X^t_t}$,
$$
\sup_{|\xi|\le R}\E[|\mun_{X^{N}_t}(e^{i\xi x} )-
\mu_{t}^\tau (e^{i\xi x})|]\le ({R\over N}+
\sup_{|\xi|\le R}\E[|\mun_{tX^N_1+(1-t)X^N_0}(e^{i\xi x} )-\nu^\tau_t(e^{i\xi x})
|]) e^{4R^2t}.$$
Therefore, if we define the distance $d_F$
on $\Pa(\R)$  by
$$d_F(\mu,\mu')=\int 
|\mu(e^{i\xi x})-\mu'(e^{i\xi x})|
e^{-4\xi^2}d\xi$$
we have proved that 
there exists a finite constant $C$ 
such that for all $t\in [0,1]$,
$$\E[d_F(\mun_{X^{N}_t},\mu_{t}^\tau)]\le C d_F(\mun_{tX^N_1+(1-t)X^N_0},
\nu^\tau_t)+{C\over N}.$$

It is not hard to convince
ourselves that $d_F$ is a distance 
compatible with the weak topology on $\Pa(\R)$.

Observe now that on $\{d(\mun_1,\mu_1)<\d,
d(\mun_0,\mu_0)<\d\}$, $(\mun_{tX_1+(1-t)X_0},t\in [0,1])$
is tight for the usual weak topology
so that for any $\e>0$ 
we can find $\kappa>0$
so that  for any $\tau$ and $t\in [0,1]$,  $D(\tau,\mun_{0,1})<\e$ implies
$$d_F(\mun_{tX^N_1+(1-t)X^N_0},
\nu^\tau_t)<\kappa.$$

Therefore, for any $t_1,\cdots,t_n\in [0,1]$,
for any $(X^N_0,X^N_1)\in\{d(\mun_0,\mu_0)<\d; d(\mun_1,\mu_1)<\d;
D(\mun_{0,1},\tau)<\e\}$, Chebyshev inequality yields
$$\P( \max_{1\le k\le
n}d_F(\mun_{X^{N}_{t_k}},\mu_{{t_k}}^\tau)>\eta|X^N_1)
\le nC( \kappa+{1\over N})$$
with $\mu^\tau_t=\mu_{X_t}$ the distribution of $X_t=tX_1+(1-t)X_0
+\sqrt{t(1-t)}S$ when the law of $(X_0,X_1)$ is $\tau$. 
Hence for any $\eta$, when $\kappa$ (i.e $\e$) is small enough and
$N$ large enough,
$$\P( \max_{1\le k\le n}d_F(\mun_{X^{N}_{t_k}},\mu_{t_k}^\tau)<\eta|X^N_1)
>{1\over 2}.$$
Hence
\begin{eqnarray*}
\P(d(\mun_1,\mu_1)<\d;
D(\mun_{0,1},\tau)<\e)
&\le & 2 \P (d(\mun_1,\mu_1)<\d;
D(\mun_{0,1},\tau)<\e  , \max_{1\le k\le n}d_F(\mun_{X^{N}_{t_k}},\mu_{{t_k}}^\tau
)<\eta).\\
\end{eqnarray*}
We
arrive at, for  $\e$ small enough
and  any $\tau\in\Ma_{0,1}$,
$$
\limsup_{N\ra\infty} {1\over N^2}\log\P(d(\mun_1,\mu_1)<\d;
D(\mun_{0,1},\tau)<\e)
\le  \limsup_{N\ra\infty} {1\over N^2}\log\P(
 \max_{1\le k\le n}d_F(\mun_{t_k},\mu_{t_k}^{\tau})<\d).$$
Using 
the large deviation upper bound for the law of $(\mun_t, t\in [0,1])$
from \cite{GZ2}, we deduce

\begin{eqnarray*}
\limsup_{N\ra\infty} {1\over N^2}\log\P(d(\mun_1,\mu_1)<\d)
&\le&-{\beta\over 2}\min_{1\le p\le M
}\inf_{\max_{1\le k\le n} d_F(\nu_{t_k},\mu^{\tau_p}_{t_k})\le \d}
 S(\nu)\\
\end{eqnarray*}
We can now let $\e$ going to zero,
and then $\d$ going to zero, and then $n$ going to infinity,
to obtain,
since $S$ is a good rate function, that
\begin{eqnarray*}\limsup_{\d\ra 0}\limsup_{N\ra\infty} {1\over N^2}\log\P(d(\mun_1,\mu_1)<\d)
&\le&
-{\beta\over 2}\inf_{
\tau :\tau\circ X_0^{-1}=\mu_0
\atop\tau\circ X_1^{-1}=\mu_1} S(\mu^\tau).
\end{eqnarray*}
Since it was also proved in \cite{GZ2}  that
$$\liminf_{\d\ra 0}\liminf_{N\ra\infty} {1\over N^2}\log\P(d(\mun_1,\mu_1)<\d)
\ge -{\beta\over 2}\inf_{
\nu_0=\mu_0
\atop\nu_1=\mu_1} S(\nu)$$
we obtain

$$\inf_{
\nu_0=\mu_0
\atop\nu_1=\mu_1} S(\nu)=\inf_{
\tau :\tau\circ X_0^{-1}=\mu_0
\atop\tau\circ X_1^{-1}=\mu_1} S(\mu^\tau).$$

Hence, if  $\mbox{FBB}(\mu_0,\mu_1)$ is the set of laws of free Brownian
bridges
between $\mu_0$ and $\mu_1$,
i.e
$$\mbox{FBB}(\mu_0,\mu_1)=\{\mu^\tau, \tau\circ X_0^{-1}=\mu_0
,\tau\circ X_1^{-1}=\mu_1\},$$
we have seen 
that
$$\inf\{S(\nu), \nu_0=\mu_0,\nu_1=\mu_1\}
=\inf\{ S(\nu), \nu\in \mbox{FBB}(\mu_0,\mu_1)\}.$$ 

To finish the proof 
of Theorem \ref{fbb}, we need to
show that $\mbox{FBB}(\mu_0,\mu_1)$
is a closed subset of $\Ca([0,1],\Pa(\R))$
so that indeed the infimum is reached in $\mbox{FBB}(\mu_0,\mu_1)$.

Observe here that $\mu^\tau$ does depend
only partially on $\tau$ since it
only depends on $\{\nu^\tau_t,
t\in [0,1]\}$.
Noting that
$$\nu^\tau_t(x^p)=\sum_{r=0}^p
t^r \tau(P_{r,p}(X_1-X_0,X_0))$$
with $P_{r,p}(X,Y)$ the sum over all
the monomial functions with total
degree $p$ and degree $r$ in $X$,
we see that $\mu^\tau$ only
depends on the restriction of $\tau$ 
to polynomial functions $P\in \Sa=\{
P_{r,p}, 0\le r\le p<\infty\}$.
Of course, 
$$\Ma^{\Sa,C}_{0,1}=\{\tau|_{\Sa}, \tau\in\Ma_{0,1},
\tau(X^{2p}+Y^{2p})\le 2C^{2p}, \forall p\in \N \}$$
is closed for the dual topology generated by 
the polynomial functions of $\Sa$. Here $C$ denotes a common 
uniform bound on $X_0$ and $X_1$, 
and we
have
$$ \mbox{FBB}(\mu_0,\mu_1)=\{\mu^{\tau|_{\Sa}}, \tau\in\Ma_{0,1}\}
=\{\mu^\kappa, \kappa \in\Ma^{\Sa,C}_{0,1}\}.  $$
We denote, for $\kappa\in \Ma^{\Sa,C}_{0,1}$
and $t\in [0,1]$, $\nu_t^\kappa\in\Pa(\R)$ the
distribution of $tX_1+(1-t)X_0$ 
when the joint distribution of $(X_0,X_1)$ restricted 
to $\Sa$ is $\kappa$. Then, $\mu^\kappa_t=\nu_t^\kappa{\tiny\boxplus}
\s_{t(1-t)}$.
We now show that $ \mbox{FBB}(\mu_0,\mu_1)$
is a closed set of $\Ca([0,1],\Pa(\R))$,
which  insures, since $S$ is a good rate function 
on  $\Ca([0,1],\Pa(\R))$, that the infimum is achieved 
on $\mbox{FBB}(\mu_0,\mu_1)$.
Indeed, if $\mu^n$ is a sequence of
$ \mbox{FBB}(\mu_0,\mu_1)$ given by $\{\nu_t^{\kappa_n}{\boxplus}
\s_{t(1-t)},t\in [0,1]\}$,
the weak convergence of $\mu^n$ 
implies the weak convergence of $\kappa^n$.
Indeed,
 for any $p\in\N$, any $t\in [0,1]$, 
$$\mu^n_t(x^p)=\nu^{\kappa_n}_t(x^p)
+P_t( \mu^n_t(x^l)
, l\le p-1)$$
with a polynomial function $P_t$.
Hence, by induction, the convergence 
of $(\mu^n_t(x^p))_{p\in\N}$ (recall that
$\mu^n$ is supported by $[-C-2,C+2]$ for
any $n$ so
that weak convergence is equivalent 
to moments convergence) results with the convergence 
of $\left(\nu^{\kappa_n}_t(x^p))\right)_{p\in\N}$,
and again, since $(\nu^{\kappa_n}_t)_{n\in\N}$
is supported by $[-C,C]$, with
the weak convergence of $\nu^{\kappa_n}_t$
towards some probability measure $\nu_t$.
Since this convergence holds for any $t\in [0,1]$,
we can expend the moments in powers of the time variable to
conclude that $\kappa_n$ converges towards $\kappa \in\Ma^{\Sa,C}_{0,1}$.
Again by free convolution calculus, this
convergence results with the convergence 
of $\mu^n$ towards $\mu^\kappa\in\mbox{FBB}(\mu_0,\mu_1) $. Hence, 
 $\mbox{FBB}(\mu_0,\mu_1)$
is closed.

\subsection{ {Proof  of  Lemma \ref{l3}:}}\label{app2}
In \cite{GZ2} (see (2.13) and  Lemma 2.10) O. Zeitouni and I
proved 
that for any path $\nu\in\Ca^1([0,1],\Pa(\R))$,
there exists a path $\nu^{\e,\D}$ such that
$$\limsup_{\e,\D\downarrow 0} S^{0,1}(\nu^{\e,\D})= S_{\mu_0}(\nu).$$
This path was constructed as follows.
Let $P_\e$ be the Cauchy law with parameter $\e$
and set $\mu^\e=P_\e* \mu$
be the convoluted path with the Cauchy law.
Moreover,
if $0=t_1<t_2<\ldots<t_n=1$ with $t_i=(i-1)\D$,
 we set, for $t\in [t_k,t_{k+1}[$,
$$
\nu^{\e,\D}_t=\nu^{\e}_{t_k}+\frac{(t-t_k)}{\D} [\nu^\e_{t_{k+1}}-
\nu^\e_{t_k}].$$

Let us therefore consider $S^{0,1}(\nu^{\e,\D})$.
Because we took the convolution with respect to the Cauchy law,
the Hilbert transform $H\nu^{\e,\D}_t$ 
is well defined, and actually a continuously differentiable function
with respect to the time variable and an analytic function
with respect to the space variable.
Henceforth, in the supremum defining $S^{0,1}(\nu^{\e,\D})$, 
we can actually make the change of function $f(t,x)
\ra f(t,x)-\int\log|x-y|d\nu^{\e,\D}_t(y)$.
Observing that, with $\nu_i^\e=\nu_i*P_\e$
for $i\in\{0,1\}$,
$$\int_0^1\int\partial_t\left(
\int\log|x-y|^{-1}d\nu^{\e,\D}_t(y)\right)
d\nu^{\e,\D}_t(x)dt={1\over 2}\left(\Sigma(\nu_1^\e)
-\Sigma(\nu_0^\e)\right),$$
we find that

\begin{eqnarray*}
S^{0,1}(\nu^{\e,\D})&=&{1\over 2}\left(\Sigma(\nu_1^\e)
-\Sigma(\nu_0^\e)\right) 
+{1\over 2}\int_0^1\int (H\nu^{\e,\D}_t)^2 d\nu^{\e,\D}_t
dt\\
&&+\sup_{f\in\Ca^{2,1}_b([0,1]\ts\R)}
\{ \int f_1d\mu_1^\e-\int f_0 d\nu_0^\e
 -\int_0^1 \int\partial_t f_t d\nu_t^{\e,\D}
dt -{1\over 2}<f,f>^{0,1}_{\nu^{\e,\D}}\}\\
&\ge& {1\over 2}\left(\Sigma(\nu_1^\e)
-\Sigma(\nu_0^\e)\right) +{1\over 2}\int_0^1\int (H\nu^{\e,\D}_t)^2 d\mu^{\e,\D}_t
dt\\
\end{eqnarray*}
where in the last line we observed that

\begin{eqnarray*}
&&\sup_{f\in\Ca^{2,1}_b([0,1]\ts\R)}
\{ \int f_1d\nu_1^\e-\int f_0 d\nu_0^\e -\int_0^1 \int\partial_t f_t 
d\nu_t^{\e,\D}
dt -{1\over 2}<f,f>^{0,1}_{\nu^{\e,\D}}\}\\
&=&\sup_{f\in\Ca^{2,1}_b([0,1]\ts\R)}\sup_{\l\in\R}
\{ \l\int f_1d\nu_1^\e-\l\int f_0 d\nu_0^\e -\l\int_0^1
 \int\partial_t f_t 
d\nu_t^{\e,\D}
dt -{\l^2\over 2}<f,f>^{0,1}_{\nu^{\e,\D}}\}\\
&=& {1\over 2}
\sup_{f\in\Ca^{2,1}_b([0,1]\ts\R)}
\left({ ( \int f_1d\nu_1^\e-\int f_0 d\nu_0^\e
 -\int_0^1 \int\partial_t f_t d\nu_t^{\e,\D})^2
\over <f,f>^{0,1}_{\nu^{\e,\D}}}\right)\\
&\ge&0
\\
\end{eqnarray*}
Observing that
$$\int_0^1\int (H\nu^{\e,\D}_t)^2 d\mu^{\e,\D}_t
dt={\D}\sum_{k=0}^{[{1\over\D}]}
\int (H\nu^{\e}_{t_k})^2 d\nu^{\e}_{t_k}$$
converges since 
$t\ra H\nu^{\e}_t$ and $t\ra \nu^{\e}_t$
are continuous for any $\nu\in\Ca([0,1],\Pa(\R))$, we arrive at
\begin{eqnarray}
\liminf_{\D\downarrow 0}
S^{0,1}(\nu^{\e,\D})&\ge &
{1\over 2}\left(\Sigma(\nu_1^\e)
-\Sigma(\nu_0^\e)\right) +{1\over 2}\int_0^1\int (H\nu^{\e}_t)^2
 d\nu^{\e}_t
dt\label{qop}\\
\nonumber
\end{eqnarray}
Notice that for $t\in\{0,1\}$,

$$\Sigma(\nu_t^\e)=\int\log|x-y|^{-1}dP_\e*\nu_t(x)dP_\e*\nu_t(y)
={1\over 2}\int \log( (x-y)^2+\e^2)^{-1} d\nu_t(x) d\nu_t(y).$$
Hence, monotone convergence  theorem asserts that
$$\lim_{\e\downarrow 0} \Sigma(\nu_t^\e)= \Sigma(\nu_t).$$
Now, recall that for any $\rho\in L_3$, 
Tricomi \cite{tricomi} p. 169 asserts that

$$ \frac{\pi^2}{2} \rho(x)^2=
\frac{1}{2} (H\rho)^2(x) - H(\rho(H\rho))(x)\,,$$
so that
$$\int (H\rho)^2 (x) \rho(x)dx={\pi^2\over 3}
\int (\rho(x))^3 dx.$$
Since, for any $\e>0$, $\nu^\e_t$ is absolutely continuous
with respect to
Lebesgue
measure with density $\rho_t^\e \in L_3(dx)$ for almost all $t\in [0,1]$,
we conclude by  Fatou's lemma that
that
\begin{eqnarray*}
+\infty>\liminf_{\e\downarrow 0}\liminf_{\D\downarrow 0}
S^{0,1}(\nu^{\e,\D})&\ge &
{1\over 2}\left(\Sigma(\nu_1)
-\Sigma(\nu_0)\right)+{\pi^2\over 6}\int_0^1\int
\liminf_{\e\downarrow
0} (\rho^\e_t(x))^3
 dx
dt.\\
\end{eqnarray*}
Finally, it is easy to see that (\ref{qop})
implies that $\mu_t(dx)=\rho_t(x)  dx$ for almost
all $t\in (0,1)$ and then that
 $\rho^\e_t(x)$ converges towards $\rho_t$
almost surely. Hence, we have proved that

\begin{eqnarray*}
S_{\mu_0}(\nu_.)
&\ge &
{1\over 2}\left(\Sigma(\nu_1)
-\Sigma(\nu_0)\right)+{\pi^2\over 6}\int_0^1\int
 (\rho_t(x))^3
 dx
dt.\\
\end{eqnarray*}

\hfill\xx

\nn
{\bf Acknowledgments :} I am very much indebted 
towards C. Villani and O. Zeitouni whose careful reading 
of preliminary versions of the manuscript, wise remarks and
encouragments  were crucial in this research. I
am also very grateful to D. Serre and Y. Brenier
for stimulating discussions.

\def\refname{BIBLIOGRAPHY}

\end{document}